\documentclass[a4paper,12pt]{article}

\usepackage[margin=1.2in]{geometry}

\usepackage{amsmath}
\usepackage{algorithmic}
\usepackage{algorithm}
\usepackage{amsthm}
\usepackage{amsfonts}
\usepackage{comment}
\usepackage{graphicx}
\usepackage{hyperref}
\usepackage{color}
\usepackage{subcaption}
\usepackage{comment}
\usepackage{makecell}
\usepackage{array}

\newtheorem{remark} {Remark}

\newtheorem{theorem} {Theorem}
\newtheorem{lemma} {Lemma}
\newtheorem{definition} {Definition}

\def\s{{\mathbf{s}}}

\def\p{{\mathbf{p}}}

\def\x{{\mathbf{x}}}

\def\u{{\mathbf{u}}}
\def\v{{\mathbf{v}}}
\def\z{{\mathbf{z}}}
\def\w{{\mathbf{w}}}
\def\y{{\mathbf{y}}}
\def\p{{\mathbf{p}}}
\def\b{{\mathbf{b}}}
\def\X{{\mathbf{X}}}
\def\Y{{\mathbf{Y}}}
\def\A{{\mathbf{A}}}
\def\M{{\mathbf{M}}}
\def\I{{\mathbf{I}}}
\def\S{{\mathbf{S}}}
\def\B{{\mathbf{B}}}

\def\V{{\mathbf{V}}}
\def\Z{{\mathbf{Z}}}
\def\W{{\mathbf{W}}}
\def\U{{\mathbf{U}}}

\newcommand{\mI}{\mathcal{I}}

\newcommand{\mV}{\mathcal{V}}

\newcommand{\mX}{\mathcal{X}}

\newcommand{\mP}{\mathcal{P}}
\newcommand{\mF}{\mathcal{F}}

\newcommand{\mK}{\mathcal{K}}
\newcommand{\mS}{\mathcal{S}}

\newcommand{\ballnuc}{\mathcal{B}_{\Vert{\cdot}\Vert_*}}
\newcommand{\mbS}{\mathbb{S}}

\newcommand{\E}{\mathbb{E}}

\newcommand{\trace}{\mathrm{Tr}}
\newcommand{\rank}{\mathrm{rank}}
\newcommand{\reals}{\mathbb{R}}

\title{Accelerated Frank-Wolfe Algorithms: Complementarity Conditions and Sparsity}
\author{Dan Garber \\ \small{dangar@technion.ac.il}}

\date{\vspace{-8pt}\small{Faculty of Data and Decision Sciences \\ \vspace{2pt} Technion - Israel Institute of Technology}}

\begin{document}

\maketitle

\begin{abstract}
We develop new accelerated first-order algorithms in the Frank-Wolfe (FW) family for minimizing smooth convex functions over compact convex sets, with a focus on two prominent constraint classes: (1) polytopes and (2) matrix domains given by the spectrahedron and nuclear-norm balls. A key technical ingredient is a complementarity condition that captures solution sparsity---face dimension for polytopes and rank for matrices. We present two algorithms: (1) a purely linear optimization oracle (LOO) method for polytopes that has optimal worst-case first-order (FO) oracle complexity and, aside of a finite \emph{burn-in} phase and up to a logarithmic factor, has LOO complexity that scales with $r/\sqrt{\epsilon}$, where $\epsilon$ is the target accuracy and $r$ is the solution sparsity (independently of the ambient dimension), and (2) a hybrid scheme that combines FW with a sparse projection oracle (e.g., low-rank SVDs for matrix domains with low-rank solutions), which also has optimal FO oracle complexity, and after a finite burn-in phase, only requires $O(1/\sqrt{\epsilon})$ sparse projections and LOO calls (independently of both the ambient dimension and the sparsity level of optimal solutions). Our results close a gap on how to accelerate recent advancements in linearly-converging FW algorithms for strongly convex optimization, without paying the price of the dimension.
\end{abstract}

\section{Introduction}
We consider algorithms based on the well known Frank-Wolfe (FW) technique \cite{frank1956algorithm, Jaggi13} for solving the constrained convex optimization problem:
\begin{align}\label{eq:genOptProb}
\min_{\x\in\mK}f(\x),
\end{align}
where $\mK\subset\E$ is a convex and compact subset of a Euclidean vector space $\E$, and $f:\E\rightarrow\reals$ is convex and $\beta$-smooth (Lipschitz-gradient) w.r.t. the Euclidean norm over $\E$. In particular, we will be interested in two important classes of constraints: (1) $\mK$ is a polytope in $\reals^n$, or (2) $\mK$ is a set of real matrices with $\ell_1$-bounded singular values, i.e., $\mK$ is either the spectrahedron in $\mbS^n$ --- the set of real symmetric positive semidefinite $n\times n$ matrices with unit trace, which will be denoted as $\mS^n$, or $\mK$ is the unit nuclear norm ball of matrices in $\reals^{m\times n}$, i.e., the set of all real $m\times n$ matrices with sum of singular values at most 1, which will be denoted as $\ballnuc^{m,n}$. We focus on these two classes of constraints for two reasons. First, these are the most prominent classes of constraints in which FW-based methods can be significantly more efficient than projection-based methods, which is indeed the reason for the surge in interest in this classical technique in recent years, as implementing the linear optimization oracle (LOO) can be much more efficient than projection, see for instance discussions in \cite{Jaggi13, COMBETTES2021565, hazan2012projection} \footnote{in the case of polytopes, this statement is not generic: for arbitrary polytopes linear optimization is not significantly more efficient than projection, however for many structured polytopes, such as those that arise from well-studied combinatorial structure, FW updates can indeed be significantly more efficient}. Second, for these types of constraints, new variants of the original FW method have been designed in recent years that are significantly faster under suitable assumptions such as strong convexity, as we further detail below.

The main drawback of the standard FW method is that it suffers from a slow worst-case sublinear convergence rate of $O(1/\epsilon)$ to reach an $\epsilon$-approximated solution (in function value) \cite{Jaggi13}. This is already sub-optimal, both in terms of number of queries to the first-order (FO) oracle of $f$ and number of optimization steps over $\mK$\footnote{we use the term optimization step loosely to refer to solving a (conceptually) simple optimization problem over $\mK$ such as linear optimization or projection}, under our assumptions on Problem \eqref{eq:genOptProb}, compared to accelerated projection-based methods which enjoy a convergence rate $O(1/\sqrt{\epsilon})$ \cite{nesterov2018lectures, beck2017first}. Moreover, if Problem \eqref{eq:genOptProb} also satisfies lower curvature conditions, such as strong convexity of $f(\cdot)$, standard projected gradient methods (either accelerated or not) converge with a linear rate, i.e., $O(\log{1/\epsilon})$, while the rate of standard FW does not improve in general under this additional strong assumption \cite{lan2013complexity}. 

For both polytopes and the trace-bounded matrix domains listed above, various works in recent years developed new algorithms based on the FW technique which are able to leverage strong convexity-like conditions to yield linear convergence rates that scale with $\log{1/\epsilon}$. For polytopes these methods all rely on the principle of introducing so-called \textit{away-steps} into the algorithm (on top of the standard updates) \cite{garber2016linearly, lacoste2015global}. While this modification indeed results in convergence rate that scales only with $\log{1/\epsilon}$, it unfortunately also scales in worst-case with the ambient dimension $n$ (which is unavoidable, e.g.,  \cite{lan2013complexity}). Some works showed that under an additional strict complementarity condition (see definition in the sequel), and provided a finite burn-in phase (which is independent of $n$), the convergence rate scales only with the dimension of the optimal face of the polytope, which can be significantly faster when this dimension is much smaller than $n$, i.e., optimal solutions are sparse, see the classical work \cite{guelat1986some} and the refined analysis in \cite{garber2020revisiting}. We also mention in passing that \cite{garber2016linear, bashiri2017decomposition} analyzed FW variants for polytopes that converge with a linear rate that scales only with the dimension of the optimal face without a burn-in phase and without assuming strict complementarity, however these hold only for highly structured polytopes named \textit{simplex-like polytopes} in \cite{bashiri2017decomposition}.

In the case of the spectrahedron or the unit nuclear norm ball, a completely different modification of the basic FW method, sometimes referred to as Block-FW methods, has been introduced that leverages strong convexity-like conditions \cite{allen2017linear, garber2021improved, ding2020spectral}. These methods replace the standard LOO, which for these matrix domains corresponds to a rank-one SVD computation (leading singular vector), with computing a rank-$r$ SVD for some positive integer $r$. Assuming all optimal solutions have rank at most $r$, this modification indeed results, under strong convexity, in a linear convergence rate (such that is independent of both the ambient dimension and $r$ and in fact matches the standard projected gradient method) \cite{allen2017linear}. Importantly, when $r$ is relatively small, these low-rank SVDs can still be far more efficient than computing the exact projection.
Here we also mention that it was established in \cite{garber2021convergence} that under a strict complementarity condition (or even weaker complementarity conditions) and at certain proximity of optimal solutions, the exact Euclidean projection itself has rank at most $r$ and thus can be computed using only a rank-$r$ SVD.

While, as surveyed above, the question of leveraging strong convexity-like conditions in FW-based methods has led to significant developments in algorithms and complexities, the question of improving the complexities without strong convexity-like assumptions and pushing it towards those of optimal accelerated projection-based methods, i.e., convergence rates (both in terms of FO queries and optimization steps) that scale only with $1/\sqrt{\epsilon}$, still presents significant challenges. One natural approach is to leverage the fact that the Euclidean projection problem itself is smooth and strongly convex, and thus we can in principle run a projection-based accelerated gradient method and use the above mentioned fast FW methods for strongly convex optimization to efficiently solve the auxiliary projection problems (to sufficient accuracy). This approach however has a severe limitation: while as surveyed above, such algorithms can leverage the sparsity of optimal solutions (whether it is sparsity in the sense of low-dimensionality of the optimal face for polytopes or low rank for matrices, and under strict complementarity in the polytope case), to obtain complexities independent of the ambient dimension,  it is not obvious that the auxiliary projection problems within accelerated methods will satisfy these sparsity conditions and retain the dimension-independent complexities.  

In this work we close this gap. We establish via customized methods that build on the advancements surveyed above and novel analysis, that the above natural approach can indeed be made to work while avoiding explicit dimension dependency, under complementarity conditions (Definition \ref{def:comp}). We design accelerated FW-based algorithms that for polytopes, the spectrahedron, and unit nuclear norm ball, require number of FO queries that scales with $O(1/\sqrt{\epsilon})$, and aside of a finite (dimension-independent) burn-in phase and up to a logarithmic factor,  require number of optimization steps that also scales with  $O(1/\sqrt{\epsilon})$. Importantly, our algorithms do not require the strongest \textit{strict complementarity} condition (e.g., as in \cite{guelat1986some, garber2020revisiting}), but are adaptive to a hierarchy of such conditions, and automatically optimize the balance between the complexity of the first slow burn-in phase and that of the second more efficient phase.  

More concretely, we present two algorithms, both build on the celebrated FISTA method \cite{beck2009fast} with a custom analysis inspired by the conditional gradient sliding framework \cite{lan2016conditional}. 
One algorithm is purely LOO-based and is intended only for polytopes. It uses the \textit{away-step Frank-Wolfe} method (AFW) \cite{lacoste2015global} (though with a different analysis) to solve auxiliary problems. Our second algorithm, while relying on a LOO, also relies on the availability of a \textit{sparse projection oracle} (Definition \ref{def:sp}). For the spectrahedron and unit nuclear norm ball this oracle amounts to computing low-rank SVDs in the same way as in Block-FW methods \cite{allen2017linear}. This algorithm is also beneficial for some polytopes, e.g., the unit simplex or $\ell_1$-balls, in which the sparse projection amounts to computing the  Euclidean projection only w.r.t. some of the top entries in the vector to project (see more details in the sequel). 

For both algorithms, we combine an Approximated-FISTA  (AFISTA) method, which provides convergence guarantees for both the function-value error and the distance between consecutive iterates, with novel identification arguments tailored to the moving auxiliary problems arising within AFISTA. This allows us to leverage the complementarity conditions and show that, after a finite number of iterations, the auxiliary problems can be solved with efficiency depending only on the sparsity level associated with the complementarity condition (in the case of strict complementarity this means simply the sparsity of optimal solutions --- dimension of the optimal face for polytopes and rank of optimal solutions for the matrix domains). This is similar to the classic linear convergence analysis of AFW for polytopes studied in \cite{guelat1986some} and later refined in \cite{garber2020revisiting}, only that while these analyses are w.r.t. the global objective function $f$ and do not yield accelerated rates, and also require the strongest strict complementarity condition, here we develop similar ideas for the auxiliary problems within AFISTA which in turn yield accelerated rates and, as mentioned above, are adaptive to a hierarchy of complementarity conditions. In case the complementarity conditions do not lead to improved complexity results, our rates in terms of FO calls and LOO calls match those of the conditional gradient sliding (CGS) method \cite{lan2016conditional}, up to a logarithmic factor in the LOO complexity (whether removing this log factor is possible or not remains an open question). 

Table \ref{table:res} compares our algorithms with CGS and the standard FW method. For instance, for a polytope with integer vertices and an algebraic description including only linear equalities and entry-wise nonnegativity constraints (so the geometric constant $\mu$ is $1$), when
$r\ll n$ and $\delta$ is large enough (relative to the target accuracy so that the
burn-in term is lower order), the LOO complexity of our purely-LOO method becomes
$\tilde{O}\!(
rD_{\mF}^2\sqrt{\beta D^2/\epsilon}
)$,
 compared with $O\!(\beta D^2/\epsilon)$ for CGS.  Likewise, for the sparse-projection
method, when the burn-in term is lower order,
the complexity becomes $O\!(\sqrt{\beta D^2/\epsilon})$
LOO calls and
$O\!(\sqrt{\beta D^2/\epsilon})$
sparse projections, again compared with $O(\beta D^2/\epsilon)$ LOO calls for CGS.

Finally, we mention that the recent works \cite{diakonikolas2020locally, carderera2021parameter}, focusing only on the case that $\mK$ is a polytope and the objective function $f$ is strongly convex, have shown that by using an away-step-based FW method for polytopes, after a finite number of iterations, Problem \eqref{eq:genOptProb} can be reduced to minimizing $f$ only over the convex-hull of a relatively small  subset of the vertices of $\mK$, and hence they can simply apply an optimal accelerated gradient method over this convex-hull. However, the length of their initial finite phase scales inversely with a parameter called the \emph{critical radius} \cite{diakonikolas2020locally} or \emph{strong Wolfe gap} \cite{carderera2021parameter}, which can be arbitrary small even for problems that are otherwise well-conditioned, e.g., optimal solutions lie in the interior of $\mK$ and far from the boundary. While the dependence on these parameters in \cite{diakonikolas2020locally, carderera2021parameter} is logarithmic, this is because they only considered the strongly convex setting. Without strong convexity this dependence becomes polynomial.\\

The rest of the paper is organized as follows. In Section \ref{sec:prel} we give notation and basic definitions. In particular we review complementarity conditions for Problem \eqref{eq:genOptProb} and the concept of sparse projections. In Section \ref{sec:FISTA} we present our Approximated-FISTA method which underlies our algorithms and present two key lemmas that connect between the complementarity conditions and the sparsity of solutions to the auxiliary optimization problems within our Approximated-FISTA scheme.
In Section \ref{sec:polyAlg} we present our first, purely LOO-based, algorithm for polytopes only and analyze its convergence guarantees. In Section \ref{sec:spProjAlg} we present our second algorithm which assumes access to both a LOO and a sparse projection oracle and analyze its convergence guarantees. Finally, in Section \ref{sec:exp} we provide a numerical illustration of the theoretical results.

\begin{table*}\renewcommand{\arraystretch}{1.7}
{\footnotesize
\begin{center}

\newcolumntype{C}[1]{>{\centering\arraybackslash}m{#1}}
\begin{tabular}{ | C{9.1em} | C{4.6em} | C{14.9em} | C{6.0em}  |} 
  \hline
  Algorithm &\#FO calls & \#LOO  calls & \#sparse proj.\\
  \hline
  Frank-Wolfe \cite{Jaggi13} &$\frac{\beta{}D^2}{\epsilon}$  & $\frac{\beta{}D^2}{\epsilon}$ &0  \\ \hline
  Conditional Gradient Sliding \cite{lan2016conditional} &$\sqrt{\frac{\beta{}D^2}{\epsilon}}$ & $\frac{\beta{}D^2}{\epsilon}$ & 0  \\ \hline
  FISTA+AFW (polytopes only) Theorem \ref{thm:main:poly} & $\sqrt{\frac{\beta{}D^2}{\epsilon}}$ & $\min\Big\{\min\{\frac{\beta^2D^2}{\delta^2}, ~\frac{\beta{}D^3}{\delta}\mu^2n\} + \sqrt{\frac{\beta{}D^2}{\epsilon}}\mu_{\mF}^2D_{\mF}^2r, ~\frac{\beta{}D^2}{\epsilon}\Big\}\log\frac{\beta{}D^2}{\epsilon}$ & 0  \\ \hline
  FISTA + FW + Sparse Proj.  Theorem \ref{thm:main:sp}& $\sqrt{\frac{\beta{}D^2}{\epsilon}}$ & $\min\Big\{\frac{\beta^2D^2}{\delta^2}\log\frac{\beta{}D^2}{\epsilon}+ \sqrt{\frac{\beta{}D^2}{\epsilon}}, ~ \frac{\beta{}D^2}{\epsilon}\log\frac{\beta{}D^2}{\epsilon}\Big\}$ & $\sqrt{\frac{\beta{}D^2}{\epsilon}}$  \\ \hline
\end{tabular}\caption{Comparison of Frank-Wolfe methods. $n$ and $D$ denote the ambient dimension and diameter of $\mK$, respectively. In case $\mK$ is a polytope, $\mu$ is a geometric constant (see \eqref{eq:mu}, e.g., for a polytope with integer vertices  and an algebraic description including only linear equalities and entry-wise nonnegativity constraints,  $\mu=1$). For a face $\mF$ of the polytope, $D_{\mF}, \mu_{\mF}$ are the diameter and geometric constant of the face ($D_{\mF} \leq D, \mu_{\mF} \leq \mu$). Parameters $r,\delta$ refer to  sparsity level and complementarity measure corresponding to a complementarity condition, respectively, see Definition \ref{def:comp}. In particular, our results hold for all suitable complementarity conditions simultaneously (see Theorems  \ref{thm:main:poly}, \ref{thm:main:sp} for exact details). The table is simplified by omitting universal constants and edge cases.}\label{table:res}
\end{center}}
\vskip -0.2in
\end{table*}\renewcommand{\arraystretch}{1.5}

\section{Notation and Preliminaries}\label{sec:prel}
\subsection{Notation}
We use lower-case boldface letters to denote vectors in some Euclidean space, e.g., $\x$, upper-case boldface letters to denote matrices in $\reals^{m\times n}$, e.g., $\A$, and lightface letters to denote scalars, e.g., $\alpha, a$. For any Euclidean space $\E$, we let $\Vert{\cdot}\Vert$ denote the Euclidean norm and $\langle{\cdot,\cdot}\rangle$ denote the standard inner-product. For matrix $\A\in\reals^{m\times n}$ we let $\Vert{\A}\Vert_2$ denote the spectral norm (largest singular value), and $\sigma_i(\A)$ denote the $i$th largest singular value. For a symmetric matrix $\B\in\mbS^n$ we let $\lambda_i(\B)$ denote the $i$th largest (signed) eigenvalue. We let $\mX^*\subseteq\mK$ denote the set of optimal solutions to Problem \eqref{eq:genOptProb} and we let $f^*\in\reals, \nabla{}f^*\in \E$ denote the corresponding optimal value and gradient direction (recall that since $f$ is smooth the gradient is constant over the optimal set $\mX^*$). We denote the Euclidean diameter of $\mK$ by $D$.

In the following let $\mP\subset\reals^n$  be a polytope in the form $\mP= \{\x\in\reals^n~|~\A_1\x=\b_1, \A_2\x\leq \b_2\}$, $\A_1\in\reals^{m_1\times n}$, $\A_2\in\reals^{m_2\times n}$, with set of vertices $\mV_{\mP}$ (i.e., $\mP =\textrm{conv}(\mV_{\mP})$, where $\textrm{conv}(\cdot)$ denotes the convex-hull).
Without loss of generality, we assume that
$\{\x\in\reals^n | \A_1\x=\b_1\}=\operatorname{aff}(\mP)$ ($\operatorname{aff}(\cdot)$ denotes the affine hull).

We recall that a face $\mF$ of $\mP$ corresponds to a subset of inequality constraints $\mI_{\mF} \subseteq [m_2]$ such that
\[
\operatorname{aff}(\mF)
=
\{\x\in\reals^n ~|~
\A_1\x=\b_1, 
\A_2(i)^\top\x=\b_2(i)~ \forall i\in\mI_{\mF}\}.
\]
The dimension of $\mF$ is given by:
\begin{align}\label{eq:faceDim}
\dim\mF
&:=
\dim\operatorname{aff}(\mF) \nonumber \\
&~=
n-\dim\operatorname{span}\left(
\{\A_1(1),\ldots,\A_1(m_1)\}
\cup
\{\A_2(i):i\in\mI_{\mF}\}
\right).
\end{align}
When taking $\mI_{\mF}=\emptyset$ we simply have
\begin{align}\label{eq:polyDim}
\dim\mP := n - \dim\textrm{span}\left(\{\A_1(1),\cdots\A_1(m_1)\}\right).
\end{align}

In case the set $\mK$ in Problem \eqref{eq:genOptProb} is a polytope we shall denote by $\mF^*$ the lowest-dimension face containing all optimal solutions. 

Following \cite{garber2016linearly}, we define a geometric constant of the polytope $\mu_{\mP}$, as follows:  we let $\mathbb{A}(\mP)$ denote the set of all $\rank(\A_2)\times n$ matrices whose rows are linearly independent rows chosen from the rows of $\A_2$, we define 
\begin{align}
\psi_{\mP} := \max_{\M\in\mathbb{A}(\mP)}\Vert{\M}\Vert_2 ,~~ \xi_{\mP} := \min_{\v\in\mV_{\mP}}\min_i\{{\b_2(i) - \A_2(i)}^{\top}\v ~|~ \b_2(i) > {\A_2(i)}^{\top}\v\},
\end{align}
where for any matrix $\A$ we let $\A(i)$ denote the column vector corresponding to the $i$th row.
We  define 
\begin{align}\label{eq:mu}
\mu_{\mP} = \left\{ \begin{array}{ll}
         \psi_{\mP}/\xi_{\mP} & \mbox{$\dim\mP >0 $};\\
        1 & \mbox{$\dim\mP = 0$}.\end{array} \right.
\end{align} 
As an example, when the vertices of $\mP$ are integer and $\A_2 =-\I, \b_2=\mathbf{0}$ (i.e., entrywise nonnegativity constraints), we simply have $\mu_{\mP} = 1$. This holds for instance for the simplex, the $[0,1]^n$ cube, the unit flow polytope, the perfect matchings polytope in a bipartite graph, and more.
In case the polytope is the set $\mK$ in Problem \eqref{eq:genOptProb} we shall simply write $\mu$. 

We denote by $D_{\mP}$ the Euclidean diameter of $\mP$.  We denote $\log_+(s) = \max\{1, \log(s)\}$.

\subsection{Complementarity conditions and sparse projections}


As detailed in the Introduction, central to our approach will be the assumption of some complementarity conditions. Such conditions are classic in the optimization literature and have been studied also in the specific context of Frank-Wolfe methods, e.g., \cite{guelat1986some, garber2020revisiting, garber2023linear, ding2020spectral, ding2020k, garber2025linearly}.  

In order to keep the presentation clear and short, we directly present for each type of constraints the complementarity conditions in the form most appropriate. For more detailed derivations of these conditions and related discussions we refer the interested reader to \cite{ding2020k}.
\begin{definition}[Complementarity conditions for polytopes, the spectrahedron, and the unit nuclear norm ball]\label{def:comp}
We shall say Problem \eqref{eq:genOptProb} satisfies the complementarity condition with dimension $r$ and complementarity measure $\delta >0$ under one of the following three cases:
\begin{itemize}
\item
$\mK$ is a polytope with a set of extreme points $\mV$, and there exists a face $\mF$ of dimension $r$ such that $\mX^*\subseteq\mF$, and
\begin{align}\label{eq:compCond:poly}
\forall \v\in\mV\setminus\mF: \quad \max_{\x^*\in\mX^*}\frac{\langle{\v-\x^*,\nabla{}f^*}\rangle}{\Vert{\v-\x^*}\Vert} \geq \delta.
\end{align}
\item
$\mK$ is the spectrahedron $\mS^n$ and 
\begin{align}\label{eq:compCond:mat}
\lambda_{n-r}(\nabla{}f^*) - \lambda_{n}(\nabla{}f^*) \geq \delta.
\end{align}
\item
$\mK$ is the unit nuclear norm ball $\ballnuc$ and
\begin{align}\label{eq:compCond:matNuc}
\sigma_{1}(\nabla{}f^*) - \sigma_{r+1}(\nabla{}f^*) \geq \delta.
\end{align}
\end{itemize}
\end{definition}
\begin{remark}
In the polytope case we shall say strict complementarity holds if \eqref{eq:compCond:poly} holds (with $\delta >0$) for the optimal face $\mF^*$ and $\mX^*\subset\mathrm{relint}(\mF^*)$. For the spectrahedron and unit nuclear norm ball we shall say strict complementarity holds if condition \eqref{eq:compCond:mat} or \eqref{eq:compCond:matNuc} hold, respectively, for some $r$ and $\delta >0$ such that any optimal solution $\x^*$ satisfies $\rank(\x^*)=r$.  We shall not require the assumption of strict complementarity in this work.
\end{remark}
\begin{remark}
For polytopes, the (strict) complementarity condition is often presented without the normalization by the distance as in \eqref{eq:compCond:poly} \cite{guelat1986some, garber2020revisiting}. However, the formulation in \eqref{eq:compCond:poly} will be more suitable for our derivations. 
\end{remark}
\begin{remark}
The unit trace/radius normalization assumed for the spectrahedron and matrix nuclear norm ball is without loss of generality: under the rescaling $\x=R\y$, the parameters transform as
$\beta\mapsto R^2\beta$, $D\mapsto D/R$, and $\delta\mapsto R\delta$, so the
quantities $\beta D^2$, $\beta D/\delta$, and $\delta^2/\beta$ appearing in
our guarantees are invariant.
\end{remark}
While some works only consider the extreme case in which \textit{strict complementarity} holds, e.g., \cite{guelat1986some, garber2020revisiting, ding2020spectral}, the conditions above are more general and allow a natural tradeoff between the dimension (or sparsity) parameter $r$ and the complementarity measure $\delta$: increasing $r$ will naturally increase the amount of computation per outer iteration of our methods, but will result in shorter burn-in time until the methods reach their "highly-efficient phase" of the run. Reducing $r$ (as long as the corresponding complementarity measure $\delta$ remains strictly positive) will naturally have the opposite effect.
Importantly, none of our algorithms will require knowledge of some complementarity measure $\delta$. Our second algorithm which is based on sparse projections (defined next) will require a target sparsity parameter $\hat{r}$ and will automatically adapt to any complementarity condition with dimension $r \leq \hat{r}$.  

\begin{definition}[sparse projection]\label{def:sp}
Given a sparsity measure $\mathrm{sp}:\mK\rightarrow\mathbb{N}_0$ and sparsity value $r\in\textrm{range}(\textrm{sp})$, we shall define the sparse projection operator $\widehat{\Pi}_{\mK}^r[\cdot]$ as
\begin{align}
\widehat{\Pi}_{\mK}^r[\x] \in\arg\min\nolimits_{\y\in\mK: \textrm{sp}(\y) \leq r}\Vert{\y-\x}\Vert.
\end{align}
Concretely, for polytopes we let $\textrm{sp}(\y)$ be the dimension of the smallest face containing $\y$ (see \eqref{eq:faceDim}), and for the spectrahedron and unit nuclear norm ball we let $\textrm{sp}(\y)=\rank(\y)$.
\end{definition}
For the spectrahedron $\mS^n$, computing $\widehat{\Pi}_{\mK}^r[\x]$ amounts to projecting (exactly) onto $\mS^n$ only the top (signed) $r$ components in the eigen-decomposition of $\x$, which in turn requires only a rank-$r$ eigen-decomposition of $\x$, which can be far more efficient than projection, which in worst-case requires full-rank eigen-decomposition, whenever $r<<n$. Similarly, for the unit nuclear norm ball $\ballnuc^{m,n}$, computing $\widehat{\Pi}_{\mK}^r[\x]$ amounts to projecting only the top $r$ components in the SVD of $\x$, see \cite{allen2017linear, garber2021improved, ding2020spectral}.
For the unit simplex polytope in $\reals^n$, computing $\widehat{\Pi}_{\mK}^r[\x]$ amounts to projecting (exactly) only the $r+1$ largest (signed) entries in $\x$ onto the unit simplex (setting other $n-r -1$ entries to zero) \cite{beck2016minimization} \footnote{Note that per \eqref{eq:faceDim}, a vector in the simplex with sparsity $s$ corresponds to a face of dimension $s-1$}.

\section{Approximated FISTA and Two Key Lemmas}\label{sec:FISTA}
As we mentioned, our algorithms build on the celebrated FISTA method \cite{beck2009fast} which we now quickly review. We consider a slightly generalized version, which to the best of our knowledge is due to \cite{chambolle2015convergence}, which will be important for our derivations. 
\begin{definition}[FISTA algorithm]\label{def:FISTA}
Fix $a\geq 2$ and let $\x_0=\y_0\in\mK$. Denote the real-valued sequence $(\lambda_t)_{t\geq 1}$ where $\lambda_t: = \frac{t+a-1}{a}$.
The FISTA algorithm produces a sequence $\{\x_t\}_{t\geq 0}$ according to the following updates:
\begin{align}
\forall t\geq 1: \quad \x_t &\gets \arg\min_{\x\in\mK}\left\{\phi_t(\x): = \langle{\x-\y_{t-1},\nabla{}f(\y_{t-1})}\rangle+\frac{\beta}{2}\Vert{\x-\y_{t-1}}\Vert^2\right\} \label{eq:FISTA:xup} \\
\y_t &\gets \x_t + \frac{\lambda_t - 1}{\lambda_{t+1}}(\x_t - \x_{t-1}) \label{eq:FISTA:yup}.
\end{align}
\end{definition}

Based on the above we consider the following Approximated-FISTA (AFISTA) scheme which allows for errors in sub-problems. Naturally, many such inexact accelerated schemes were proposed before, e.g., \cite{NIPS2011_8f7d807e}, however the one considered here, for which we do not claim particular novelty,  is carefully tailored to our needs.
\begin{definition}[Approximated-FISTA (AFISTA)]\label{def:AppFISTA}
Let $(\nu_t)_{t\geq 1}\subset\reals_+$ be a sequence of error tolerances and for any iteration $t\geq 1$ define the function 
\begin{align}\label{eq:omega}
\omega_t(\x) := \max_{\w\in\mK}\left\langle{\x - [(1-\lambda_{t}^{-1})\x_{t-1} + \lambda_{t}^{-1}\w], \nabla\phi_t(\x)}\right\rangle,
\end{align}
where $\phi_t(\cdot), \lambda_t$ are as in Definition \ref{def:FISTA}.

The AFISTA algorithm is the same as FISTA with the modification that the exact computation in  \eqref{eq:FISTA:xup} is replaced with the condition:
\begin{align}\label{eq:AFISTA:cond}
 \x_t&\in\{\x\in\mK ~|~ \omega_t(\x) \leq \nu_t\}.
 \end{align}

\end{definition}

Let $\Pi_{\mK}[\cdot]$ denote the exact Euclidean projection operator onto $\mK$. For any iteration $t\geq 1$ of AFISTA we define the following quantities which will be central throughout the rest of this work:
\begin{align}
h_t &:= f(\x_t) - f^*, \\
\x_t^* &:= \arg\min\nolimits_{\x\in\mK}\phi_t(\x) = \Pi_{\mK}\left[{\y_{t-1}-\beta^{-1}\nabla{}f(\y_{t-1})}\right], \\
d_t^* &:= \frac{1}{2}\Vert{\x_t^* - \x_{t-1}}\Vert^2, \\
d_t &:= \frac{1}{2}\Vert{\x_t - \x_{t-1}}\Vert^2.
\end{align}

The proof of the following theorem mostly adapts the analysis from \cite{chambolle2015convergence} to also account for the approximation errors in AFISTA. In particular, aside from the standard convergence rate w.r.t. function values, this theorem also bounds the sequences of distances $d_t, d_t^*$ which will be important for our derivations. 
\begin{theorem}\label{thm:AFISTA}
Consider Algorithm AFISTA with $a=5$ and let $D_0 \geq \min_{\x^*\in\mX^*}\Vert{\x_0-\x^*}\Vert$. Then, 
\begin{align*}
\forall t\geq 2:\quad h_{t} \leq \frac{\beta{}D_0^2}{2\lambda_{t}^2} +  \frac{1}{\lambda_t^2}\sum_{\tau=1}^t\lambda_{\tau}^2\nu_{\tau}, ~~ \max\{d_{t}^*, d_t\} &\leq  \frac{D_0^2}{\lambda_{t}^2} + \frac{3}{\beta\lambda_t^2}\sum_{\tau=1}^{t}\lambda_{\tau}^2\nu_{\tau}.
\end{align*}
In particular, fixing $T\geq 2$ and setting 
\begin{align}
\nu_t = \frac{\beta{}D_0^2}{\lambda_t^2t(1+\log{T})} \quad  \forall t\leq T
\end{align}
gives,
\begin{align*}
\forall 2\leq t\leq T: \quad h_t \leq \frac{3\beta{}D_0^2}{2\lambda_{t}^2}, ~~ \max\{d_t^*, d_t\} \leq \frac{4D_0^2}{\lambda_{t}^2}.
\end{align*}
\end{theorem}
\begin{proof}
The bound on $h_t$ follows from Lemma \ref{lem:FISTA:funcConv} and the bounds on $d_t^*$ and $d_t$ follow from Lemma \ref{lem:FISTA:distConv}.
\end{proof}

As mentioned in the Introduction, in order to get our complexity results, and in particular the complexity results for the burn-in phase (which, up to a logarithmic factor, is independent of the target accuracy $\epsilon$), we combine the AFISTA scheme with a technique inspired by the conditional gradient sliding method \cite{lan2016conditional}. Towards this, for any iteration $t\geq 1$ of AFISTA we define the function:
\begin{align}\label{eq:Phi_t}
\Phi_t(\w) := \lambda_t^{-1}\langle{\w-\y_{t-1}, \nabla{}f(\y_{t-1})}\rangle + \frac{\lambda_t^{-2}\beta}{2}\Vert{\w + \lambda_t\left({(1-\lambda_t^{-1})\x_{t-1}-\y_{t-1}}\right)}\Vert^2.
\end{align}
The following very simple lemma shows that approximately minimizing $\Phi_t$ over $\mK$ is sufficient to guarantee the approximation condition \eqref{eq:AFISTA:cond} in AFISTA.  The benefit of minimizing $\Phi_t$ (as opposed to directly working with the function $\phi_t$ defined in \eqref{eq:FISTA:xup}) is that as the number of iteration $t$ increases, $\Phi_t$ becomes more and more smooth (recall $\lambda_t = \Theta(t)$), and so it becomes more and more efficient to optimize with FW.
\begin{lemma}\label{lem:slide}
Fix iteration $t\geq 1$ of AFISTA and let $\w_t\in\mK$ be such that
\begin{align*}
\max_{\w\in\mK}\langle{\w_t - \w, \nabla\Phi_{t}(\w_t)}\rangle \leq \nu_t,
\end{align*}
and define $\x_t = (1-\lambda_t^{-1})\x_{t-1}  + \lambda_t^{-1}\w_t$. 
Then, $\omega_t(\x_t) \leq \nu_t$.
\end{lemma}
\begin{proof}
\begin{align*}
\omega_t(\x_t) &= \max_{\w\in\mK}\langle{\x_t - \left({(1-\lambda_t^{-1})\x_{t-1} + \lambda_t^{-1}\w}\right), \nabla\phi_t(\x_t)}\rangle \\
&=\max_{\w\in\mK}\langle{\left({(1-\lambda_t^{-1})\x_{t-1} + \lambda_t^{-1}\w_t}\right) - \left({(1-\lambda_t^{-1})\x_{t-1} + \lambda_t^{-1}\w}\right), \nabla\phi_t(\x_t)}\rangle  \\
&=\max_{\w\in\mK}\langle{\w_t - \w, \lambda_t^{-1}\nabla\phi_t(\x_t)}\rangle  \\
&= \max_{\w\in\mK}\langle{\w_t-\w,  \lambda_t^{-1}\nabla{}f(\y_{t-1}) + \beta \lambda_t^{-1}(\x_t - \y_{t-1})}\rangle \\
&= \max_{\w\in\mK}\langle{\w_t-\w,  \lambda_t^{-1}\nabla{}f(\y_{t-1}) + \beta \lambda_t^{-1}\left({(1-\lambda_t^{-1})\x_{t-1}+\lambda_t^{-1}\w_t - \y_{t-1}}\right)}\rangle\\
&= \max_{\w\in\mK}\langle{\w_t-\w,  \lambda_t^{-1}\nabla{}f(\y_{t-1}) + \beta \lambda_t^{-2}\left({\w_t + \lambda_t\left({(1-\lambda_t^{-1})\x_{t-1} - \y_{t-1}}\right)}\right)}\rangle\\
&= \max_{\w\in\mK}\langle{\w_t-\w,  \nabla\Phi_t(\w_t)}\rangle \leq \nu_t.
\end{align*}
\end{proof}

\subsection{Two key lemmas: sparsity in auxiliary problems}
The following two key lemmas, one for the case that $\mK$ is a polytope and the other for the case that it is the spectrahedron or the unit nuclear norm ball, will enable us to argue that, under a $(r,\delta)$ complementarity condition (Definition \ref{def:comp}), after a finite number of iterations, which scales inversely with $\delta$, our algorithms will automatically adapt to the sparsity level $r$. For polytopes this means, that all auxiliary problems in AFISTA will be automatically solved w.r.t. a face of dimension at most $r$, and for the matrix domains this will mean that the $r$-sparse projection will in fact be the exact Euclidean projection.

\begin{lemma}\label{lem:poly:goodSteps}
Suppose $\mK$ is a polytope in $\reals^n$ and suppose the complementarity condition \eqref{eq:compCond:poly} holds with some parameters $r,\delta$, and let $\mF$ be the corresponding face of $\mK$.
There exists a universal constant $c >0$ such that for any iteration $t\geq 2$ for which 
\begin{align}\label{eq:goodStepsCond}
\max\{h_{t-1}, \beta{}d_{t-1}, \beta{}d_t^*\} \leq \frac{c\delta^2}{\beta{}},
\end{align}
we have that $\x_t^*\in\mF$. If additionally, $\lambda_{t}^{-1} < \frac{\delta}{4\beta{}D}$, we also have that
\begin{align}\label{eq:lem:poly:goodStep:FW}
\forall \w\in\mK:\quad \arg\min\nolimits_{\u\in\mK}\langle{\u, \nabla\Phi_t(\w)}\rangle \subseteq\mF.
\end{align}
\eqref{eq:lem:poly:goodStep:FW} further implies (via the first-order optimality condition) that $\arg\min_{\w\in\mK}\Phi_t(\w) \in \mF$.
\end{lemma}
\begin{proof}
Recall that
\begin{align*}
\nabla\phi_t(\x) = \nabla{}f(\y_{t-1}) + \beta(\x- \y_{t-1}).
\end{align*}
Note that using Eq. \eqref{eq:FISTA:yup} we have that,
\begin{align*}
\Vert{\y_{t-1} - \x_{t-1}}\Vert = \frac{\lambda_{t-1}-1}{\lambda_{t}}\Vert{\x_{t-1} - \x_{t-2}}\Vert \leq \sqrt{2d_{t-1}}.
\end{align*}
Thus,
\begin{align}\label{eq:lem:scPoly:1}
\Vert{\nabla\phi_t(\x_t^*) - \nabla{}f^*}\Vert &\leq \Vert{\nabla{}f(\y_{t-1}) - \nabla{}f^*}\Vert  + \beta\Vert{\x_t^*-\y_{t-1}}\Vert \nonumber \\
&\leq \Vert{\nabla{}f(\y_{t-1}) - \nabla{}f^*}\Vert  + \beta\left({\Vert{\x_t^*-\x_{t-1}}\Vert + \Vert{\y_{t-1}-\x_{t-1}}\Vert}\right) \nonumber \\
&\leq \Vert{\nabla{}f(\y_{t-1}) - \nabla{}f^*}\Vert + \beta\left({\sqrt{2d_t^*} + \sqrt{2d_{t-1}}}\right) \nonumber \\
&\leq \Vert{\nabla{}f(\x_{t-1}) - \nabla{}f^*}\Vert + \beta\left({\sqrt{2d_t^*} + 2\sqrt{2d_{t-1}}}\right) \nonumber \\
&\leq \sqrt{2\beta{}h_{t-1}}+ \beta\left({\sqrt{2d_t^*} + 2\sqrt{2d_{t-1}}}\right),
\end{align}
where the last inequality uses a well known result for smooth and convex functions, see for instance Lemma \ref{lem:smoothIneq} in the appendix.

Fix some $\v\in\mV\setminus\mF$ and $\x_{\v}^*\in\arg\max_{\x^*\in\mX^*}\frac{\langle{\v-\x^*,\nabla{}f^*}\rangle}{\Vert{\v-\x^*}\Vert}$. Using the complementarity condition \eqref{eq:compCond:poly}, we have that
\begin{align*}
\langle{\v - \x_{\v}^*,\nabla\phi_t(\x_t^*)}\rangle &\geq \langle{\v - \x_{\v}^*,\nabla{}f^*}\rangle - \Vert{\v-\x_{\v}^*}\Vert\Vert{\nabla\phi_t(\x_t^*) - \nabla{}f^*}\Vert \\
& \geq \delta\Vert{\v-\x_{\v}^*}\Vert- \Vert{\v-\x_{\v}^*}\Vert\left({\sqrt{2\beta{}h_{t-1}}+ \beta\left({\sqrt{2d_t^*} + 2\sqrt{2d_{t-1}}}\right)}\right).
\end{align*}
Thus, under Condition \eqref{eq:goodStepsCond} we have that 
\begin{align*}
\forall \v\in\mV\setminus\mF: \quad \langle{\v - \x_{\v}^*,\nabla\phi_t(\x_t^*)}\rangle > 0.
\end{align*}
This implies that $\x_t^*\in\mF$, since otherwise, if $\x_t^*$ is supported on some vertex $\v\in\mV\setminus\mF$, the above inequality implies that the value $\phi_t(\x_t^*)$ could be further decreased (by considering the feasible point $\x_t^* + \gamma(\x_{\v}^*-\v)$ for sufficiently small positive $\gamma$), contradicting the optimality of $\x_t^*$.

We continue to prove \eqref{eq:lem:poly:goodStep:FW}. Fixing some $\eta\in(0,1]$, let us define the function 
\begin{align*}
\phi_{t,\eta}(\w) = \eta\langle{\w-\y_{t-1},\nabla{}f(\y_{t-1})}\rangle + \frac{\beta\eta^2}{2}\Vert{\w + \eta^{-1}\left({(1-\eta)\x_{t-1} - \y_{t-1}}\right)}\Vert^2,
\end{align*}
and note that $\Phi_t(\cdot) \equiv \phi_{t,\lambda_{t}^{-1}}(\cdot)$ (i.e., setting $\eta = \lambda_t^{-1}$). 

The gradient of $\nabla\phi_{t,\eta}(\w)$ is given by
\begin{align*}
\nabla\phi_{t,\eta}(\w) = \eta\nabla{}f(\y_{t-1}) + \eta^2\beta\left({\w + \eta^{-1}\left({(1-\eta)\x_{t-1} - \y_{t-1}}\right)}\right),
\end{align*}
and so,
\begin{align*}
\Vert{\nabla\phi_{t,\eta}(\w) - \eta\nabla{}f^*}\Vert &\leq  \eta\Vert{\nabla{}f(\y_{t-1})- \nabla{}f^*}\Vert \\
&~~~+ \eta\beta\Vert{\x_{t-1}-\y_{t-1}}\Vert + \eta^2\beta\Vert{\w-\x_{t-1}}\Vert \\
&\leq  \eta\Vert{\nabla{}f(\x_{t-1})- \nabla{}f^*}\Vert \\
&~~~+ 2\eta\beta\Vert{\x_{t-1}-\y_{t-1}}\Vert + \eta^2\beta\Vert{\w-\x_{t-1}}\Vert \\
&\leq  \eta\sqrt{2\beta{}h_{t-1}} + 2\eta\beta\sqrt{2d_{t-1}} + \eta^2\beta{}D,
\end{align*}
where in the last inequality we used Lemma \ref{lem:smoothIneq} again.

Fix again some $\v\in\mV\setminus\mF$ and let $\x^*_{\v}$ be as defined above. Using  again the complementarity condition \eqref{eq:compCond:poly}, we have that
\begin{align*}
\min_{\u\in\mK}\langle{\u-\v, \nabla\phi_{t,\eta}(\w)}\rangle &\leq \langle{\x_{\v}^*-\v, \nabla\phi_{t,\eta}(\w)}\rangle \\
&\leq  \langle{\x_{\v}^*-\v, \eta\nabla{}f^*}\rangle + \Vert{\v-\x_{\v}^*}\Vert\Vert{\nabla\phi_{t,\eta}(\w) - \eta\nabla{}f^*}\Vert \\
&\leq -\eta\delta\Vert{\v-\x_{\v}^*}\Vert + \Vert{\v-\x_{\v}^*}\Vert\Vert{\nabla\phi_{t,\eta}(\w) - \eta\nabla{}f^*}\Vert \\
&\leq -\eta\Vert{\v-\x_{\v}^*}\Vert\left({\delta - \sqrt{2\beta{}h_{t-1}} - 2\beta\sqrt{2d_{t-1}} - \eta\beta{}D}\right).
\end{align*}
Thus, setting $\eta = \lambda_t^{-1}$, we have that under Condition \eqref{eq:goodStepsCond} and the assumption on $\lambda_t^{-1}$, it holds that
\begin{align*}
\forall \v\in\mV\setminus\mF: \quad \min_{\u\in\mK}\langle{\u-\v, \nabla\Phi_t(\w)}\rangle < 0,
\end{align*}
which proves \eqref{eq:lem:poly:goodStep:FW}.
\end{proof}

\begin{lemma}\label{lem:goodProj}
Suppose $\mK$ is either the spectrahedron $\mS^n$ or the unit nuclear norm ball $\ballnuc^{m,n}$, and that either complementarity condition \eqref{eq:compCond:mat} or complementarity condition \eqref{eq:compCond:matNuc} holds with parameters $r,\delta$.
For any iteration $t\geq 2$ for which Condition \eqref{eq:goodStepsCond} holds, we have that $\rank(\x_t^*) \leq r$.
\end{lemma}
\begin{proof}
We prove for the spectrahedron, the proof for the unit nuclear norm ball follows the same lines with the obvious changes.

Starting from Eq.  \eqref{eq:lem:scPoly:1} in the proof of Lemma \ref{lem:poly:goodSteps} (note the derivation of \eqref{eq:lem:scPoly:1} is generic and did not rely on the specific structure of the feasible set $\mK$), we have that under the assumption of the lemma
\begin{align*}
\Vert{\nabla\phi_t(\x_t^*) - \nabla{}f^*}\Vert < \delta /2.
\end{align*}
This implies via Weyl's inequality for the eigenvalues  that
\begin{align}\label{eq:eigsGap}
&\lambda_{n-r}(\nabla\phi_t(\x_t^*)) - \lambda_{n}(\nabla\phi_t(\x_t^*)) \geq \nonumber \\
&\lambda_{n-r}(\nabla{}f^*) - \lambda_{n}(\nabla{}f^*) - 2\Vert{\nabla\phi_t(\x_t^*) - \nabla{}f^*}\Vert > 0 .
\end{align}
The first-order optimality condition implies that all eigenvectors of $\x_t^*$ with nonzero eigenvalue must also be eigenvectors of $\nabla\phi_t(\x_t^*)$ corresponding to the smallest eigenvalue $\lambda_{n}(\nabla\phi_t(\x_t^*))$. Hence, the spectral gap \eqref{eq:eigsGap} indeed implies that $\rank(\x_t^*) \leq r$, see also Lemma 5.2 in \cite{garber2021convergence}.


For the unit nuclear norm ball the proof differs only by applying Weyl's inequality for the singular values $\sigma_1$ and $\sigma_{r+1}$ (instead of eigenvalues $\lambda_{n-r}$, $\lambda_n$), and using the first-order optimality condition to establish that all singular vector pairs of $\x_t^*$ (with nonzero singular value) must be also leading singular vector pairs of $-\nabla\phi_t(\x_t^*)$, see also Lemma 2.2 in  \cite{garber2021convergence}. 
\end{proof}

\begin{remark}[tightness of Condition \eqref{eq:goodStepsCond}]
The dependence on $\delta^2/\beta$ in Condition
\eqref{eq:goodStepsCond} is unavoidable, up to universal constants.
Indeed, fix any $D,\beta>0$ and $0<\delta\leq\beta D/2$, and consider
the one-dimensional problem
\[
\mK=[0,D], \qquad f(x)=\delta x.
\]
The function $f$ is convex and $\beta$-smooth, its unique optimal
solution is $x^*=0$, and the complementarity condition
\eqref{eq:compCond:poly} holds for the optimal face $\mF^*=\{0\}$ with
complementarity measure $\delta$, since
\[
\frac{\langle D-x^*,\nabla f^*\rangle}{|D-x^*|}
=\delta.
\]

Consider iteration $t=2$ of AFISTA with $x_0=x_1=\frac{2\delta}{\beta}$. This is a valid AFISTA iterate: since $\lambda_1=1$ and $y_0=x_0$, we have
$\omega_1(x_1)
=\max_{w\in[0,D]}\delta(x_1-w)
=2\delta^2/\beta$,
and hence $x_1=x_0$ satisfies \eqref{eq:AFISTA:cond} by taking
$\nu_1=2\delta^2/\beta$.
Since $\lambda_1=1$, we have
$y_1=x_1$, and so
\[
x_2^*
=
\arg\min\nolimits_{x\in[0,D]}
\left\{
\delta(x-y_1)+\frac{\beta}{2}(x-y_1)^2
\right\}
=
\frac{\delta}{\beta}>0.
\]
Thus, $x_2^*\notin\mF^*$. However, a simple calculation gives 
$h_1=\frac{2\delta^2}{\beta}$,
 $d_1=0$, 
$\beta d_2^*=\frac{\delta^2}{2\beta}$, 
and therefore
\[
\max\{h_1,\beta d_1,\beta d_2^*\}
=
\frac{2\delta^2}{\beta}.
\]
Thus, even with the universal constant $2$ in the right-hand side of
Condition \eqref{eq:goodStepsCond}, the identification conclusion need
not hold.
\end{remark}

\section{LOO-based Algorithm for Polytopes}\label{sec:polyAlg}
In this section we present our first algorithm which is a purely LOO-based algorithm for polytopes only. The algorithm simply applies the Away-Step Frank-Wolfe method \cite{lacoste2015global} to solve the sub-problems within AFISTA on each iteration $t$, by minimizing $\Phi_t(\w)$ over $\mK$.  One crucial modification is in the initialization: we begin with a vertex minimizing the inner product with $\nabla\Phi_t(\x_{t-1})$. This is important, so when the conditions of Lemma \ref{lem:poly:goodSteps} hold, the algorithm effectively operates only on a restricted face of the polytope (the one corresponding to the complementarity condition in Lemma \ref{lem:poly:goodSteps}).
For completeness the algorithm is brought below as Algorithm \ref{alg:awayfw}.

\begin{algorithm}
\caption{Away-Step Frank-Wolfe for Polytopes}
\label{alg:awayfw}
\begin{algorithmic}[1]
\STATE input: $\x_{t-1}, \y_{t-1}, \lambda_t$, error-tolerance $\nu_t$
\STATE $\w_1 \gets \arg\min_{\u\in\mV}\langle{\u, \nabla\Phi_t(\x_{t-1})}\rangle$
\FOR{$i=1,2\dots $}
\STATE let $\w_i = \sum_{j=1}^m\rho_j\v_j$ be a convex decomposition of $\w_i$ to vertices in $\mV$, i.e., $\{\v_1,\dots,\v_m\}\subseteq\mV$, $(\rho_1,\dots,\rho_m)$ is in the unit simplex and $\forall j\in[m]:\rho_j>0$ \COMMENT{maintained explicitly throughout the run of the algorithm by tracking the vertices that enter and leave the decomposition}
\STATE $\u_i \gets \arg\min_{\v\in\mV}\langle{\v, \nabla\Phi_t(\w_i)}\rangle$, $j_i \gets \arg\max_{j\in[m]}\langle{\v_j, \nabla{}\Phi_t(\w_i)}\rangle$, $\z_i \gets \v_{j_i}$
\IF{$\langle{\w_i-\u_i, \nabla{}\Phi_t(\w_i)\rangle} \leq  \nu_t$} 
\RETURN $\x_t = (1-\lambda_{t}^{-1})\x_{t-1} + \lambda_t^{-1}\w_i$
\ENDIF
\IF{$\langle{\u_i-\w_i, \nabla{}\Phi_t(\w_i)}\rangle < \langle{\w_i-\z_i, \nabla\Phi_t(\w_i)}\rangle$} 
\STATE $\s_i \gets \u_i-\w_i$, $\gamma_{\max} \gets 1$  \COMMENT{Frank-Wolfe direction}
\ELSE
\STATE $\s_i \gets \w_i - \z_i$, $\gamma_{\max}\gets \rho_{j_i}/(1-\rho_{j_i})$ \COMMENT{away direction}
\ENDIF
\STATE $\w_{i+1} \gets \w_i+\gamma_i\s_i$ where $\gamma_i \gets\arg\min_{\gamma\in[0,\gamma_{\max}]}\Phi_t(\w_i + \gamma\s_i)$
\ENDFOR
\end{algorithmic}
\end{algorithm} 

The following theorem gives complexity bounds for Algorithm \ref{alg:awayfw}. A proof is given in the appendix. While the linear convergence rate of Algorithm \ref{alg:awayfw} was established in \cite{lacoste2015global}, here we provide a somewhat different analysis which does not rely on the \textit{pyramidal width} quantity, which is often difficult to evaluate, but follows the analysis in \cite{garber2016linearly}. Moreover, we also establish a new dimension-independent dual convergence result (the first term inside the $\min$ in \eqref{eq:afw:conv}) which is crucial to our analysis.
\begin{theorem}\label{thm:afw}
Assume $\mK$ is a polytope in $\reals^n$ and fix some iteration $t$ of AFISTA. Algorithm \ref{alg:awayfw} stops after at most
\begin{align}\label{eq:afw:conv}
O\left({\min\left\{\frac{\beta{}D^2}{\lambda_t^2\nu_t}, ~\max\left\{1, \mu^2D^2\dim\mK\right\}\log_+\left({\frac{\beta{}D^2}{\lambda_t^2\nu_t}}\right)\right\}}\right)
\end{align}
iterations, and the returned point $\x_t$ satisfies $\omega_t(\x_t) \leq \nu_t$. 

Moreover, if for some face $\mF$ of $\mK$ it holds that $\arg\min_{\w\in\mK}\Phi_t(\w)\in\mF$, $\w_1\in\mF$ and $\u_i\in\mF$ for all $i\geq 1$, then the term $\mu^2D^2\dim\mK$ in \eqref{eq:afw:conv} could be replaced with $\mu_{\mF}^2D_{\mF}^2\dim\mF$.
\end{theorem}

We can now finally present our first main result: the convergence guarantees of AFISTA when using Algorithm \ref{alg:awayfw} for computing the feasible iterates $(\x_t)_{t\geq 1}$.
\begin{theorem}\label{thm:main:poly}
Suppose $\mK$ is a polytope in $\reals^n$ and fix $\epsilon > 0$. There exists $T_{\epsilon} = O\left({\sqrt{\beta{}D^2/\epsilon}}\right)$ such that running AFISTA  with $a=5$ and $D_0 =D$ for $T_{\epsilon}$ iterations, where on each iteration $t$, $\x_t$ is computed via Algorithm \ref{alg:awayfw} with error tolerance $\nu_t$ as prescribed in Theorem \ref{thm:AFISTA}, guarantees that $h_{T_{\epsilon}} \leq \epsilon$ and the overall number of LOO calls is 
\begin{align}\label{eq:main:poly:res1}
O\left({\min\left\{\frac{\beta{}D^2}{\epsilon}, ~\sqrt{\frac{\beta{}D^2}{\epsilon}}\max\left\{1, \mu^2D^2{}n\right\}\right\}\log\frac{\beta{}D^2}{\epsilon}}\right).
\end{align}
Moreover, assuming the complementarity condition \eqref{eq:compCond:poly} holds with some parameters $r, \delta$ and letting $\mF$ be the corresponding face, and assuming $\epsilon$ is small enough (so all assumptions of Lemma \ref{lem:poly:goodSteps} hold), we have that the overall number of LOO calls is 
\begin{align}\label{eq:main:poly:res2}
&O\left({\min\left\{\left({\frac{\beta{}D}{\delta}}\right)^2,~\frac{\beta{}D}{\delta}\max\left\{1, \mu^2D^2{}n\right\}\right\}\log\frac{\beta{}D^2}{\epsilon}}\right) \nonumber \\
&~ + O\left({\sqrt{\frac{\beta{}D^2}{\epsilon}}\max\left\{1, \mu_{\mF}^2D_{\mF}^2{}\dim\mF\right\}\log\frac{\beta{}D^2}{\epsilon}}\right).
\end{align}
\end{theorem}
Note that \eqref{eq:main:poly:res2} holds simultaneously for all faces $\mF$ satisfying the complementarity condition  \eqref{eq:compCond:poly}. Hence, for a fixed $\epsilon >0$,  the LOO-complexity of the algorithm is upper-bounded by the minimum of  \eqref{eq:main:poly:res2} over all such faces.
\begin{proof}[Proof of Theorem \ref{thm:main:poly}]
The bound on $T_{\epsilon}$ follows immediately from Theorem \ref{thm:AFISTA}. We obtain the bound in \eqref{eq:main:poly:res1} by an immediate application of Theorem \ref{thm:afw} w.r.t. the polytope $\mK$. Indeed, with this theorem, we have the overall number of calls to the LOO can be upper-bounded by:
\begin{align}\label{eq:main:poly:1}
&\sum_{t=1}^{T_{\epsilon}}O\left({\min\left\{\frac{\beta{}D^2}{\lambda_t^2\nu_t}, ~\max\left\{1, \mu^2D^2n\right\}\log\left({\frac{\beta{}D^2}{\lambda_t^2\nu_t}}\right)\right\}}\right) \nonumber \\
&= \sum_{t=1}^{T_{\epsilon}}O\left({\min\left\{t{}\log{}T_{\epsilon}, ~\max\left\{1, \mu^2D^2n\right\}\log(T_{\epsilon}\log{T_{\epsilon}})\right\}}\right),
\end{align}
which yields \eqref{eq:main:poly:res1} after plugging-in the bound on $T_{\epsilon}$ and slightly simplifying. 

To prove \eqref{eq:main:poly:res2}, we first observe that as an immediate consequence of Theorem \ref{thm:AFISTA}, we have that if $T_{\epsilon} \geq T_0 = \frac{c\beta{}D}{\delta}$, for an appropriate universal constant $c$, then for all $t \geq T_0$, the conditions of Lemma \ref{lem:poly:goodSteps} hold and thus, for all $t \geq T_0$, each invocation of Algorithm \ref{alg:awayfw} acts as if optimizing only over the face $\mF$. Thus, by applying Theorem \ref{thm:afw} w.r.t. the face $\mF$, we immediately recover the second term in the sum in  \eqref{eq:main:poly:res2}. The first term in the sum in \eqref{eq:main:poly:res2} follows from the number of LOO calls until iteration $T_0$ and is upper bounded exactly as in Eq. \eqref{eq:main:poly:1}.

\end{proof}

\section{Sparse Projections-based Algorithm}\label{sec:spProjAlg}
In this section we present our second algorithm which is suitable whenever, on top of access to a LOO for $\mK$, we also have access to an oracle computing sparse projections onto $\mK$ (Definition \ref{def:sp}), which is in particular suitable when $\mK$ is the spectrahedron $\mS^n$ or the unit nuclear norm ball $\ballnuc^{m,n}$, and Problem \eqref{eq:genOptProb} has only low-rank solutions. Our algorithm uses Algorithm \ref{alg:sparseProjFW} below to compute the sequence of feasible points $(\x_t)_{t\geq 1}$ within AFISTA.

Algorithm \ref{alg:sparseProjFW} has two steps. First, it attempts to directly minimize the auxiliary function $\phi_t$ within AFISTA using the sparse projection oracle with a pre-specified sparsity level $\hat{r}$. It then uses a single call to the LOO to validate whether the sparse projection is sufficiently accurate (by checking if $\omega_t(\x) \leq \nu_t$).  If not, it simply runs the standard FW with line-search in order to minimize $\Phi_t(\w)$ over $\mK$, which in turn leads to sufficient accuracy w.r.t $\phi_t$ (via Lemma \ref{lem:slide}).

\begin{algorithm}
\caption{Sparse projection or Frank-Wolfe}
\label{alg:sparseProjFW}
\begin{algorithmic}
\STATE input: $\x_{t-1}, \y_{t-1}, \lambda_t$, sparsity parameter $\hat{r}$, error-tolerance $\nu_t$
\STATE $\x \gets \widehat{\Pi}_{\mK}^{\hat{r}}[\y_{t-1} - \beta^{-1}\nabla{}f(\y_{t-1})]$ 
\IF{$\omega_t(\x) \leq \nu_t$}
\RETURN $\x_t = \x$
\ENDIF
\STATE $\w_1 \gets \x$ \COMMENT{sparse projection failed, start standard FW iterations}
\FOR{$i=1,2,...$} 
\STATE $\u_i \gets \arg\min_{\u\in\mK}\langle{\u, \nabla\Phi_{t}(\w_i)}\rangle$
\IF {$\langle{\w_i - \u_i, \nabla\Phi_t(\w_i)}\rangle \leq \nu_t$} 
\RETURN $\x_t= (1-\lambda_t^{-1})\x_{t-1} + \lambda_t^{-1}\w_i$ \COMMENT{by Lemma \ref{lem:slide}, $\omega_t(\x_t) \leq \nu_t$}
\ELSE
\STATE $\gamma_i \gets \arg\min_{\gamma\in[0,1]}\Phi_t((1-\gamma)\w_i + \gamma\u_i)$
\STATE $\w_{i+1} \gets (1-\gamma_i)\w_i + \gamma_i\u_i$
\ENDIF
\ENDFOR
\end{algorithmic}
\end{algorithm} 

The following theorem follows directly from the dual convergence of the standard FW method, combined with Lemma \ref{lem:slide}.
\begin{theorem}[\cite{Jaggi13}, Theorem 2]\label{thm:FW}
Fix some iteration $t$ of AFISTA and some error parameter $\nu_t$. Algorithm \ref{alg:sparseProjFW} stops after $O\left({\frac{\beta{}D^2}{\lambda_t^2\nu_t}}\right)$ iterations, and the returned point $\x_t$ satisfies $\omega_t(\x_t) \leq \nu_t$. 
\end{theorem}

We are now ready to present our second main result.
\begin{theorem}\label{thm:main:sp}
Suppose $\mK$ is either a polytope in $\reals^n$, the spectrahedron $\mS^n$ or the unit nuclear norm ball $\ballnuc^{m,n}$ and fix $\epsilon > 0$. There exists $T_{\epsilon} = O\left({\sqrt{\beta{}D^2/\epsilon}}\right)$  such that running AFISTA  with $a=5$ and $D_0 =D$ for $T_{\epsilon}$ iterations, where on each iteration $t$, $\x_t$ is computed via Algorithm \ref{alg:sparseProjFW} with some sparsity parameter $\hat{r}$ and with error tolerance $\nu_t$ as prescribed in Theorem \ref{thm:AFISTA}, guarantees that $h_{T_{\epsilon}} \leq \epsilon$ and the overall number of LOO calls is 
\begin{align}\label{eq:main:sp:res1}
O\left({\frac{\beta{}D^2}{\epsilon}\log\frac{\beta{}D^2}{\epsilon}}\right).
\end{align}
Moreover, assuming one of the complementarity conditions \eqref{eq:compCond:poly}, \eqref{eq:compCond:mat} or \eqref{eq:compCond:matNuc} holds (with compatibility with the structure of $\mK$) with parameters $(r, \delta)$ such that $r \leq \hat{r}$, and assuming $\epsilon$ is small enough (so the  assumptions of either Lemma \ref{lem:poly:goodSteps} or Lemma \ref{lem:goodProj} hold), we have that the overall number of LOO calls is only
\begin{align}\label{eq:main:sp:res2}
O\left({\left({\frac{\beta{}D}{\delta}}\right)^2\log\frac{\beta{}D^2}{\epsilon} + \sqrt{\frac{\beta{}D^2}{\epsilon}}}\right),
\end{align}
and  after $O\left({\frac{\beta{}D}{\delta}}\right)$ AFISTA iterations we always have that $\widehat{\Pi}_{\mK}^{\hat{r}}[\y_{t-1} - \beta^{-1}\nabla{}f(\y_{t-1})] = \Pi_{\mK}[\y_{t-1} - \beta^{-1}\nabla{}f(\y_{t-1})] = \arg\min_{\x\in\mK}\phi_t(\x)$ (i.e., the $\hat{r}$-sparse projection is the exact projection), and the for-loop in Algorithm \ref{alg:sparseProjFW} is no longer executed.
\end{theorem}
Note that for a fixed sparsity parameter $\hat{r}$, the guarantees of the algorithm hold simultaneously w.r.t. all complementarity conditions (\eqref{eq:compCond:poly}, \eqref{eq:compCond:mat} or \eqref{eq:compCond:matNuc}) of dimension $r \leq \hat{r}$. Hence, for a fixed $\epsilon >0$,  the LOO-complexity of the algorithm is upper-bounded by the minimum of  \eqref{eq:main:sp:res2} over all such conditions.
\begin{proof}[Proof of Theorem \ref{thm:main:sp}]
The bound on $T_{\epsilon}$ follows immediately from Theorem \ref{thm:AFISTA}. We obtain the bound in \eqref{eq:main:sp:res1} by an immediate application of Theorem \ref{thm:FW} which yields that the overall number of calls to the LOO can be upper-bounded by:
\begin{align}\label{eq:main:sp:1}
&\sum_{t=1}^{T_{\epsilon}}O\left({\frac{\beta{}D^2}{\lambda_t^2\nu_t}}\right) = \sum_{t=1}^{T_{\epsilon}}O\left({t{}\log{}T_{\epsilon}}\right)
= O\left({T_{\epsilon}^2\log{T_{\epsilon}}}\right),
\end{align}
which yields \eqref{eq:main:sp:res1} after plugging-in the bound on $T_{\epsilon}$ and simplifying. 

To prove the second part of the theorem, we first observe that as an immediate consequence of Theorem \ref{thm:AFISTA}, we have that if $T_{\epsilon} \geq T_0 = \frac{c\beta{}D}{\delta}$, for an appropriate universal $c$, then for all $t \geq T_0$, the conditions of Lemma \ref{lem:poly:goodSteps} (if $\mK$ is a polytope) or Lemma  \ref{lem:goodProj} (if $\mK$ is $\mS^n$ or $\ballnuc^{m,n}$) hold and thus, for all $t \geq T_0$, each invocation of Algorithm \ref{alg:sparseProjFW} indeed terminates after the sparse projection step without entering the for-loop (since then $\x$ in particular minimizes $\phi_t$ over $\mK$ and so, from the first-order optimality condition, $\omega_t(\x) \leq 0$), and thus each such invocation makes a single call to the LOO (to evaluate $\omega_t(\x)$), which yields the second term in the sum in \eqref{eq:main:sp:res2}.
It remains to upper-bound the number of calls to the LOO until iteration $T_0$ is reached, which corresponds to the first term in the sum in \eqref{eq:main:sp:res2}. This follows exactly as in \eqref{eq:main:sp:1} by summing only over the first $T_0$ summands, instead of all $T_{\epsilon}$. 
\end{proof}

\begin{remark}[replacing FW with AFW for polytopes]\label{rem:1}In case $\mK$ is a polytope, we may want to replace the standard FW iterations in Algorithm \ref{alg:sparseProjFW} with the AFW iterations of Algorithm \ref{alg:awayfw}, since for polytopes these are often expected to converge faster. This will not deteriorate the complexity results in Theorem \ref{thm:main:sp} since Theorem \ref{thm:FW}, which is used to bound the number of FW steps, could  readily be replaced with Theorem \ref{thm:afw}.
\end{remark}

\pagebreak 
\section{Numerical Illustration}\label{sec:exp}
The focus of this work is on the theoretical results. The experiments in this section are intended only to illustrate the qualitative behavior predicted by our analysis, rather than to provide a comprehensive empirical evaluation. We therefore consider a simple synthetic setting in which the relevant sparsity and complementarity parameters can be directly controlled.
We consider  the problem of minimizing a convex quadratic function over the unit simplex  $\Delta_n = \{\x\in\reals^n~|~\x\geq \mathbf{0}, ~\mathbf{1}^{\top}\x = 1\}$:
\begin{align}\label{eq:expProb}
\min_{\x\in\Delta_n}\{f(\x) := \frac{1}{2}\x^{\top}\A\x + \b^{\top}\x\}.
\end{align}
We first note that while projection onto the simplex is quite efficient ($O(n\log{}n)$ time), the main benefit of FW methods for Problem \eqref{eq:expProb} is that, as opposed to projection-based methods which often have dense updates and so computing the gradient direction for \eqref{eq:expProb} requires $O(n^2)$ time (with dense $\A$), FW methods (including our purely LOO-based and sparse projections-based) only update a small  number of coordinates per iteration: single coordinate for purely LOO methods and $\hat{r}$ coordinates for a sparse projection step with sparsity level $\hat{r}$, and thus the time to update the gradient of  \eqref{eq:expProb} is only $O(n)$ per FW/AFW step, and $O(n\hat{r})$ per sparse projection step.

We let $\A$ be a random symmetric positive definite matrix with largest eigenvalue $\beta$.  For a selected sparsity value $r$ we let $\x^*$ be a random $r$-sparse vector in $\Delta_n$. Finally, for a desired strict complementarity measure $\delta$, we set $\b = -\A\x^* + \delta\z^*$, where for all $i\in[n]$ we have $\z^*(i) = 0$ if $i$ is in the support of $\x^*$ and $\z^*(i)=1$ otherwise. Letting $S$ denote the support of $\x^*$, this guarantees that 
\begin{equation}\label{eq:sc}
    \bigl(\nabla f(\x^*)\bigr)_i
    \;=\; 0 \quad \text{for } i\in S,
    \qquad
    \bigl(\nabla f(\x^*)\bigr)_i
    \;=\; \delta \quad \text{for } i\notin S,
\end{equation}
which in turn implies that $\x^*$ is indeed an optimal solution which satisfies strict complementarity with measure $\delta$.

We compare our two AFISTA implementations: (1) when using  AFW  to solve the inner optimization problems (AFISTA-AFW), and (2) when using the sparse projection-based Algorithm \ref{alg:sparseProjFW}, and when per Remark \ref{rem:1} we replace the standard FW iterations in  Algorithm \ref{alg:sparseProjFW} with AFW (AFISTA-SP/AFW), to the Conditional Gradient Sliding Algorithm, implemented exactly as in \cite{lan2016conditional} (CGS), and the vanilla Frank-Wolfe with line-search algorithm. Our algorithms have been implemented exactly as stated (in particular with $a=5$ in the FISTA $(\lambda_t)_{t\geq 1}$ sequence, and the sequence $(\nu_t)_{t\geq 1}$ listed in Theorem \ref{thm:AFISTA}). For Algorithm \ref{alg:sparseProjFW} we simply set $\hat{r} =r$, i.e., we use exact knowledge of the sparsity of optimal solutions.

For all algorithms we measure the convergence rate vs. the number of outer-iterations (i.e., number of sub-problems solved, for vanilla Frank-Wolfe this is just the standard iterations). Additionally, we measure the approximation error vs. the number of calls to the LOO. Since, as discussed above, the time to  update the gradient direction is proportional to the number of LOO calls, this gives a credible implementation independent  estimate for the runtime of the algorithms. In case of AFISTA-SP/AFW which also computes $r$-sparse projections, which  produces $r$-sparse vectors (and hence the time to update the gradient is $O(nr)$), we count each such call as $r$ calls to the LOO.

We set $n = 200$, $r \in \{10, 20, 40, 80\}$, $\delta \in \{0.0, 0.1, 1.0\}$, and $\beta = 100$. We use $T=2000$ outer iterations and each plot is the average of 10 i.i.d. runs.

\begin{figure}[H]
     \centering
     \begin{subfigure}[b]{0.32\textwidth}
         \centering
         \includegraphics[width=\textwidth]{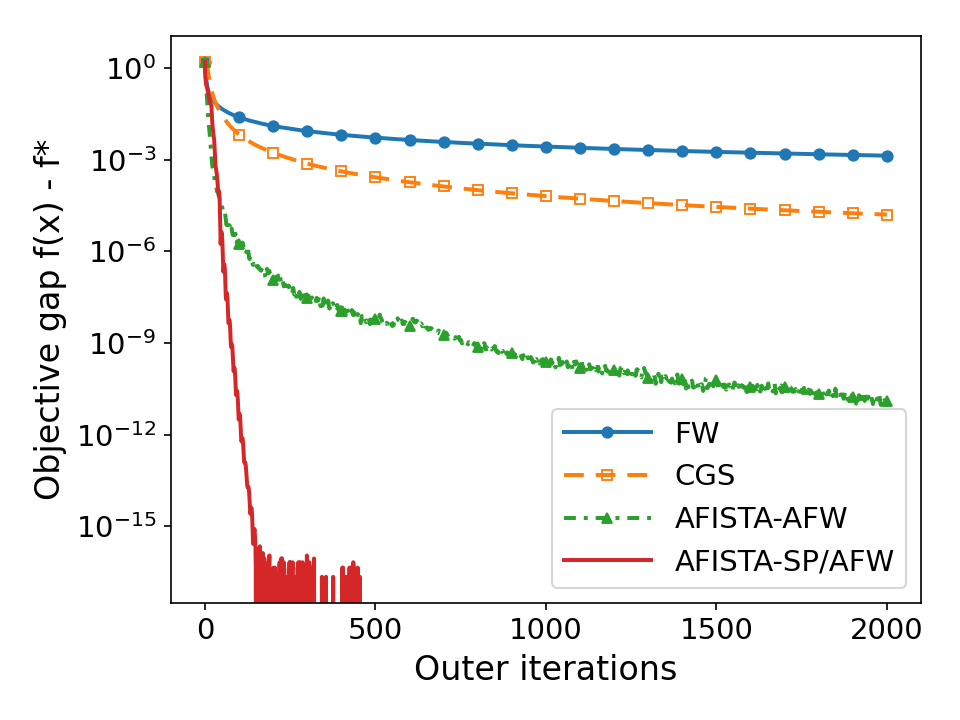}
         \caption*{$r = 10, \delta=0.0$}
         \label{fig:y equals x}
     \end{subfigure}
     \hfill
     \begin{subfigure}[b]{0.32\textwidth}
         \centering
         \includegraphics[width=\textwidth]{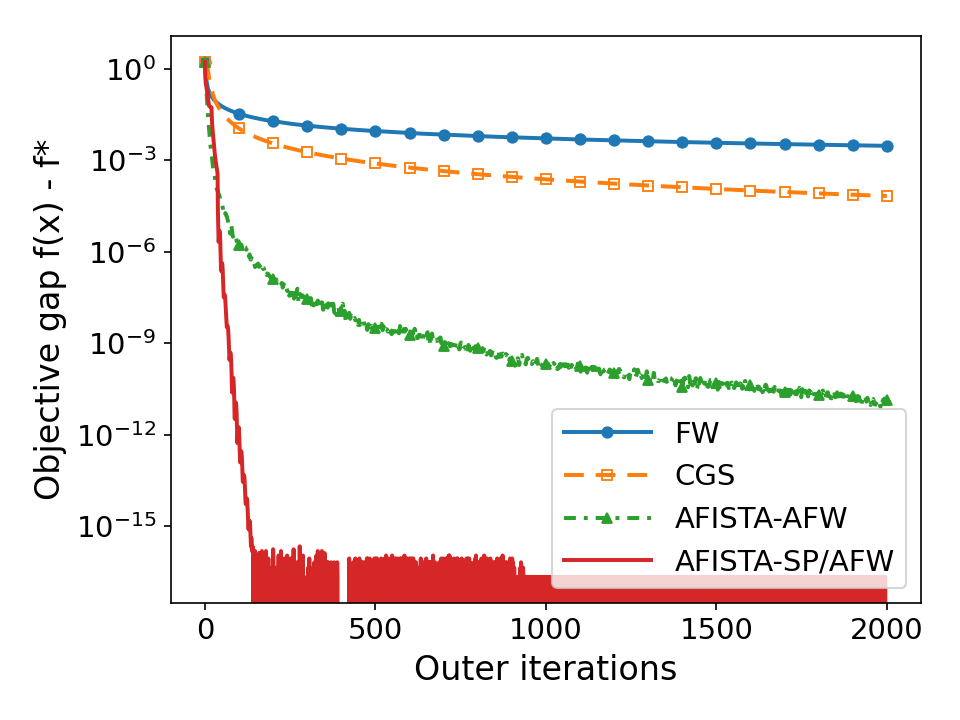}
         \caption*{$r = 10, \delta=0.1$}
         \label{fig:three sin x}
     \end{subfigure}
     \hfill
     \begin{subfigure}[b]{0.32\textwidth}
         \centering
         \includegraphics[width=\textwidth]{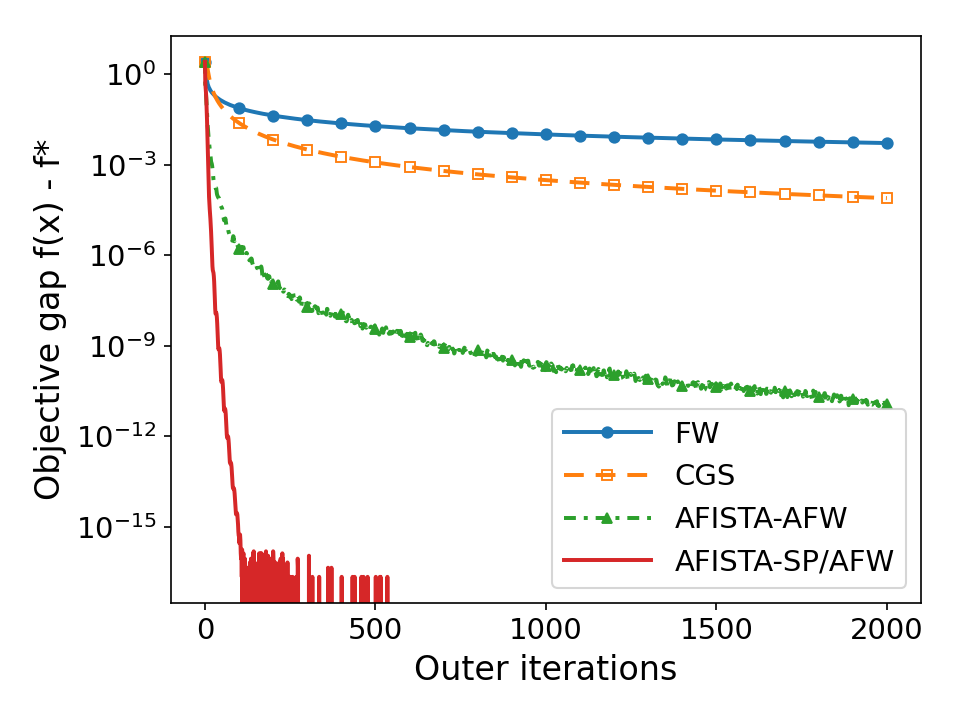}
         \caption*{$r = 10, \delta=1.0$}
         \label{fig:five over x}
     \end{subfigure}

     \begin{subfigure}[b]{0.32\textwidth}
         \centering
         \includegraphics[width=\textwidth]{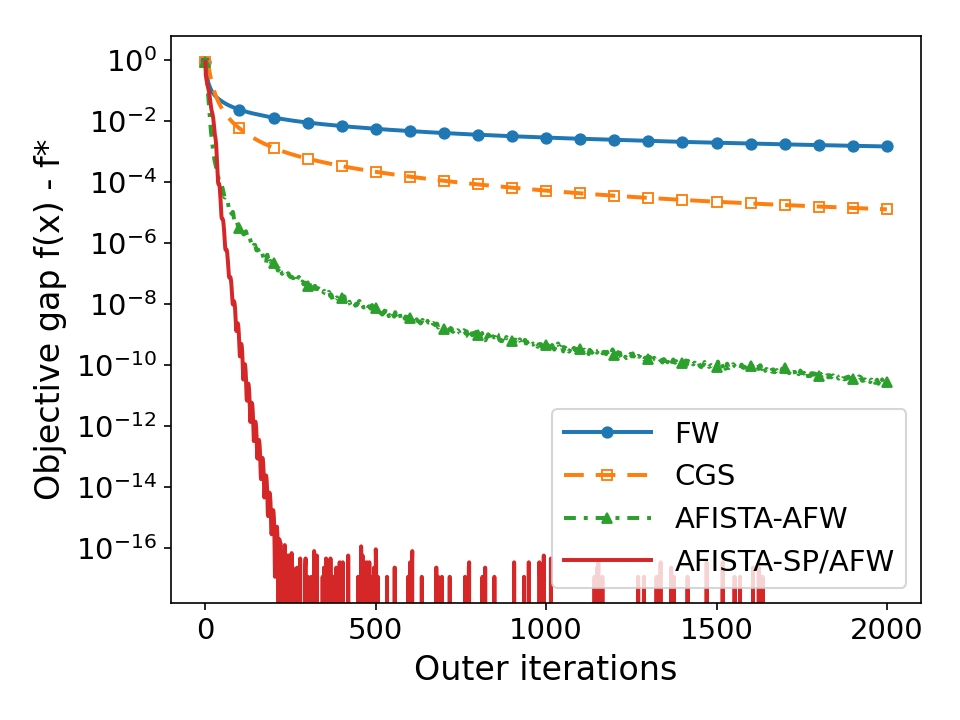}
         \caption*{$r = 20, \delta=0.0$}
         \label{fig:y equals x}
     \end{subfigure}
     \hfill
     \begin{subfigure}[b]{0.32\textwidth}
         \centering
         \includegraphics[width=\textwidth]{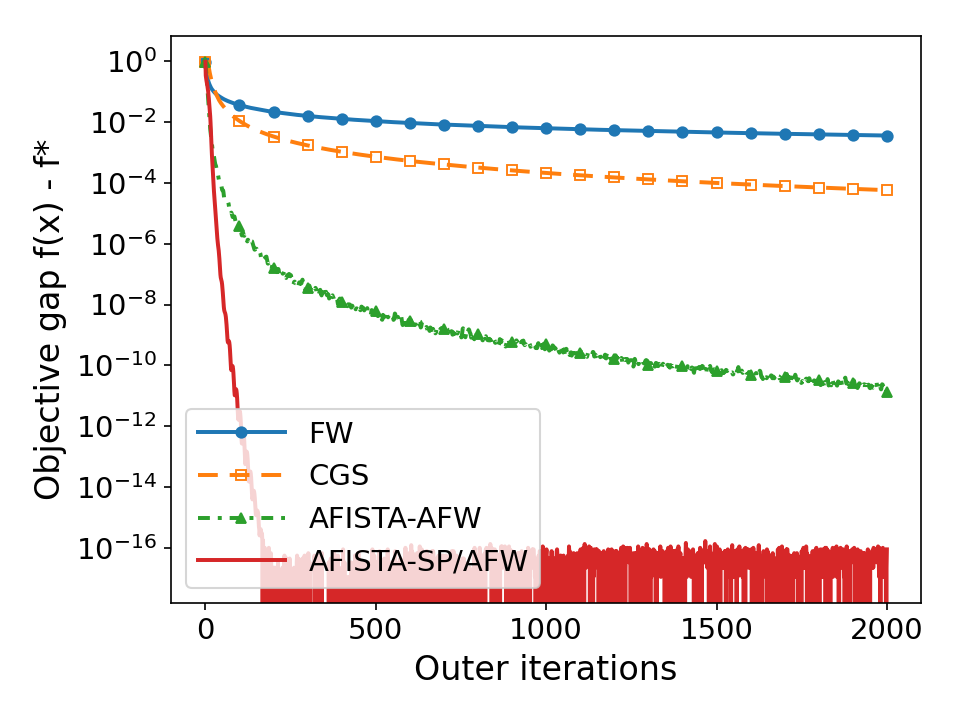}
         \caption*{$r = 20, \delta=0.1$}
         \label{fig:three sin x}
     \end{subfigure}
     \hfill
     \begin{subfigure}[b]{0.32\textwidth}
         \centering
         \includegraphics[width=\textwidth]{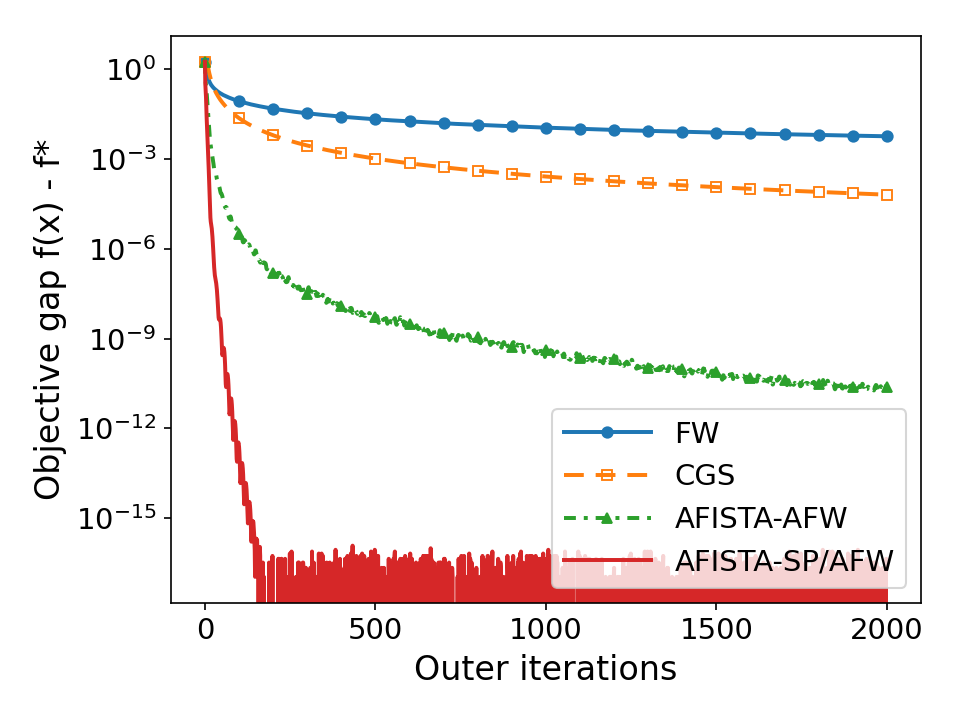}
         \caption*{$r = 20, \delta=1.0$}
         \label{fig:five over x}
     \end{subfigure}

     \begin{subfigure}[b]{0.32\textwidth}
         \centering
         \includegraphics[width=\textwidth]{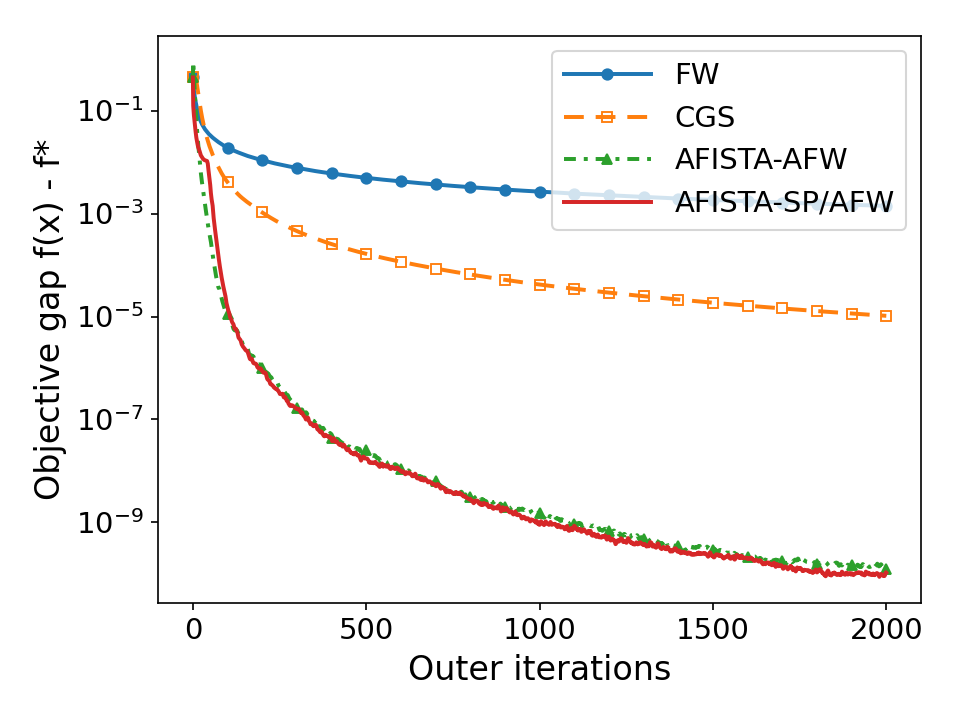}
         \caption*{$r = 40, \delta=0.0$}
         \label{fig:y equals x}
     \end{subfigure}
     \hfill
     \begin{subfigure}[b]{0.32\textwidth}
         \centering
         \includegraphics[width=\textwidth]{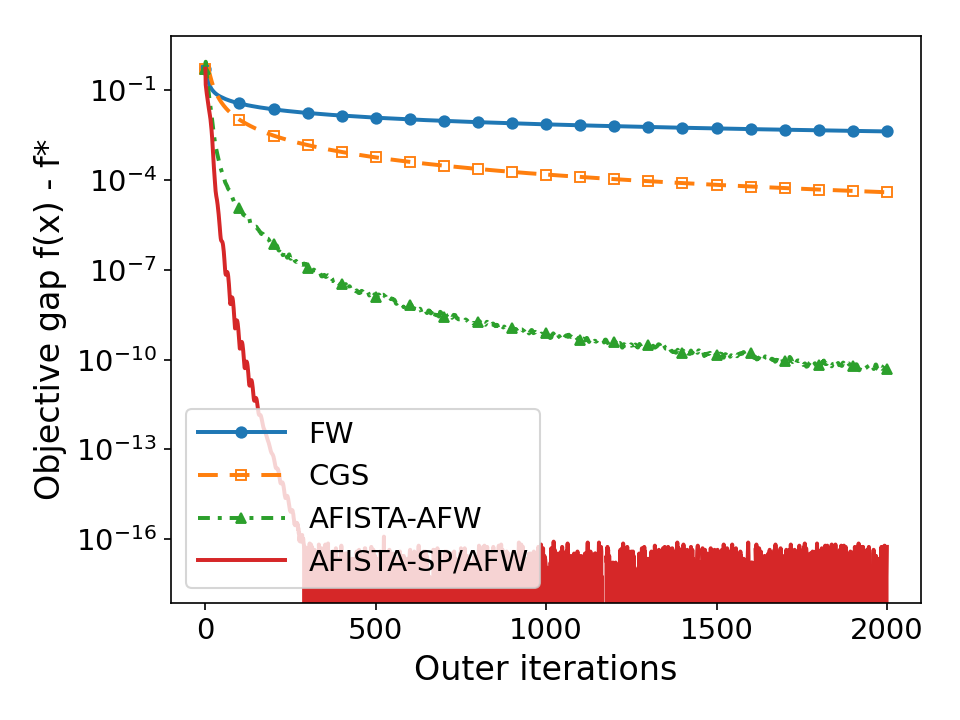}
         \caption*{$r = 40, \delta=0.1$}
         \label{fig:three sin x}
     \end{subfigure}
     \hfill
     \begin{subfigure}[b]{0.32\textwidth}
         \centering
         \includegraphics[width=\textwidth]{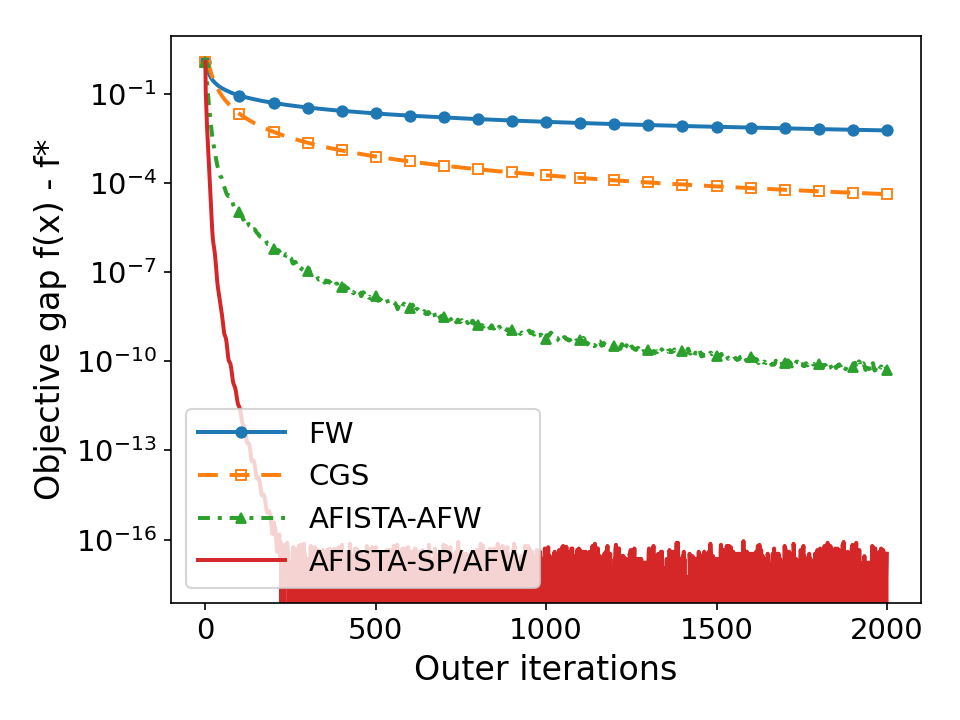}
         \caption*{$r = 40, \delta=1.0$}
         \label{fig:five over x}
     \end{subfigure}

     \begin{subfigure}[b]{0.32\textwidth}
         \centering
         \includegraphics[width=\textwidth]{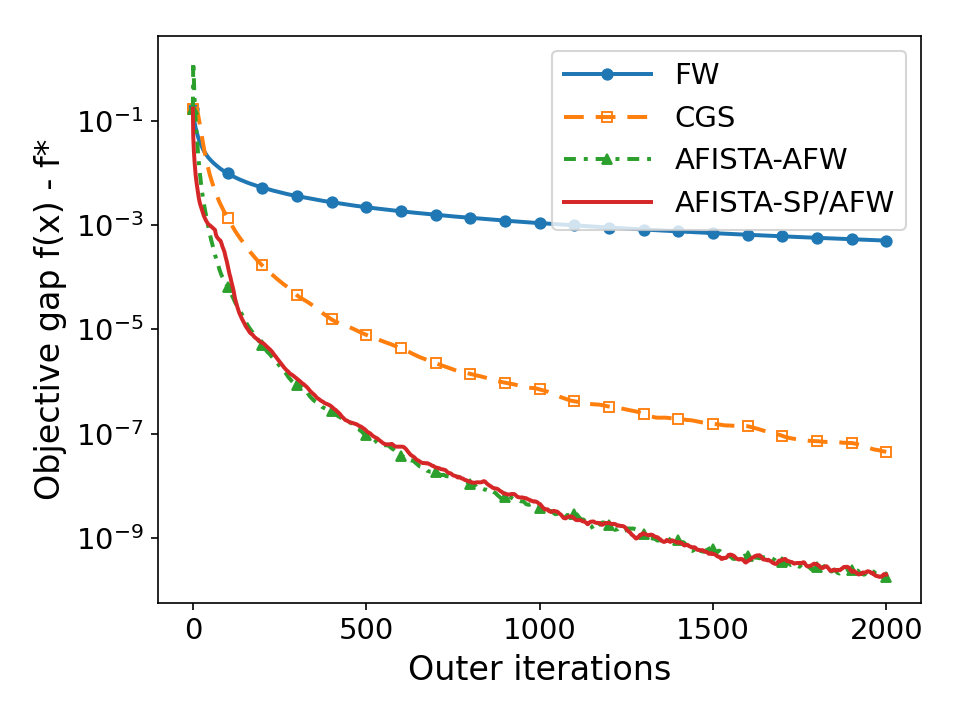}
         \caption*{$r = 80, \delta=0.0$}
         \label{fig:y equals x}
     \end{subfigure}
     \hfill
     \begin{subfigure}[b]{0.32\textwidth}
         \centering
         \includegraphics[width=\textwidth]{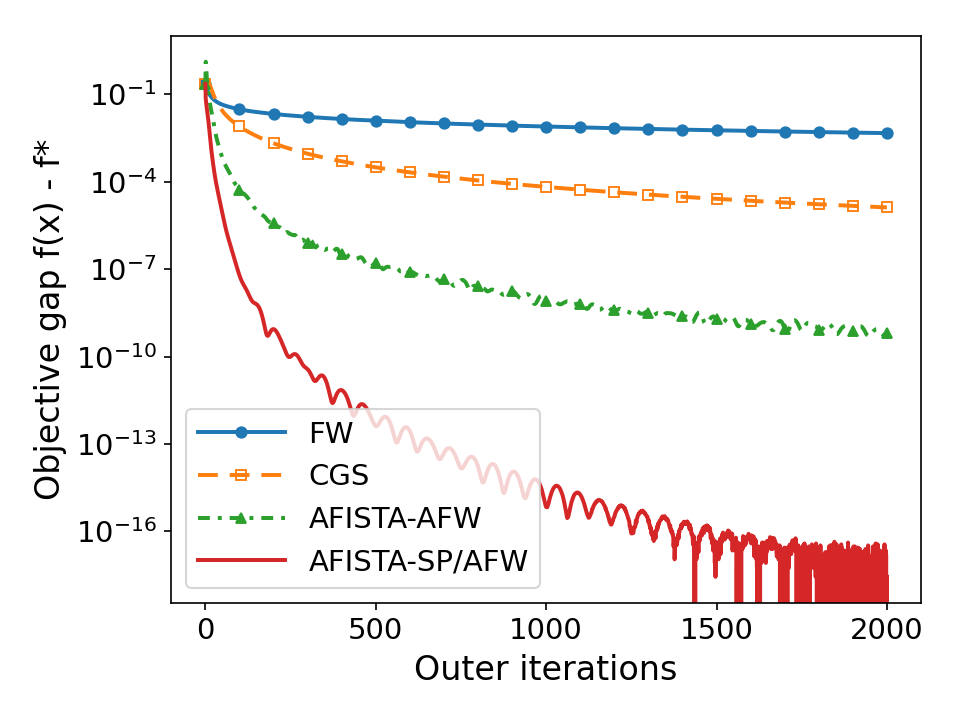}
         \caption*{$r = 80, \delta=0.1$}
         \label{fig:three sin x}
     \end{subfigure}
     \hfill
     \begin{subfigure}[b]{0.32\textwidth}
         \centering
         \includegraphics[width=\textwidth]{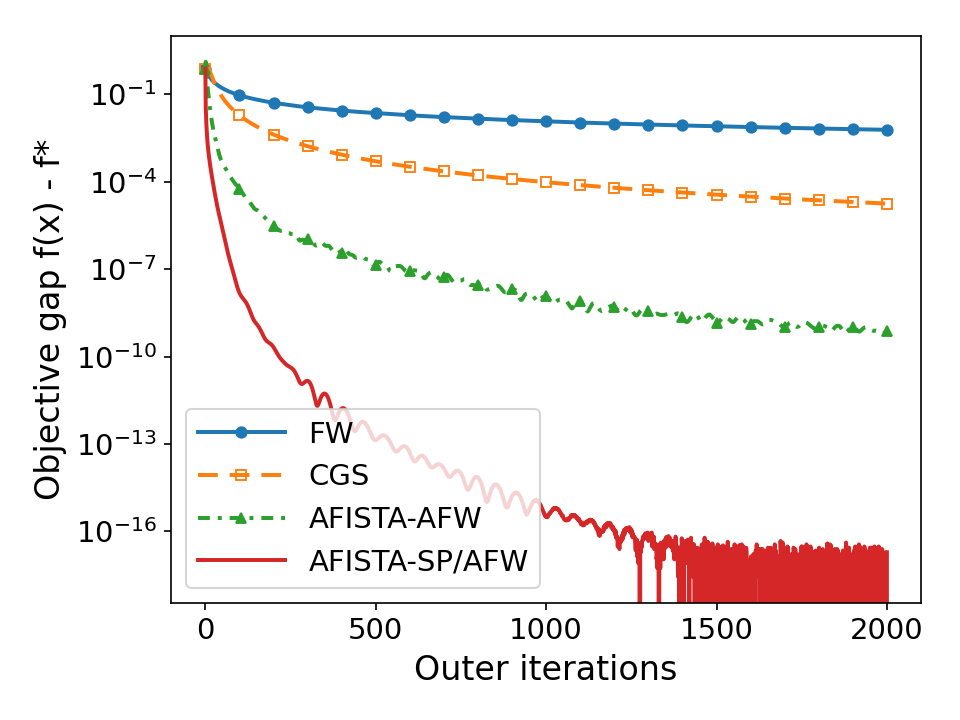}
         \caption*{$r = 80, \delta=1.0$}
         \label{fig:five over x}
     \end{subfigure}
        \caption{Approximation errors vs. number of outer iterations.}
        \label{fig:outer}
\end{figure}

We clearly see in Figure \ref{fig:outer} that both of our algorithms clearly dominate in terms of convergence w.r.t. number of outer iterations, where AFISTA-SP/AFW is most often significantly faster than all other methods.

\begin{figure}[H]
     \centering
     \begin{subfigure}[b]{0.32\textwidth}
         \centering
         \includegraphics[width=\textwidth]{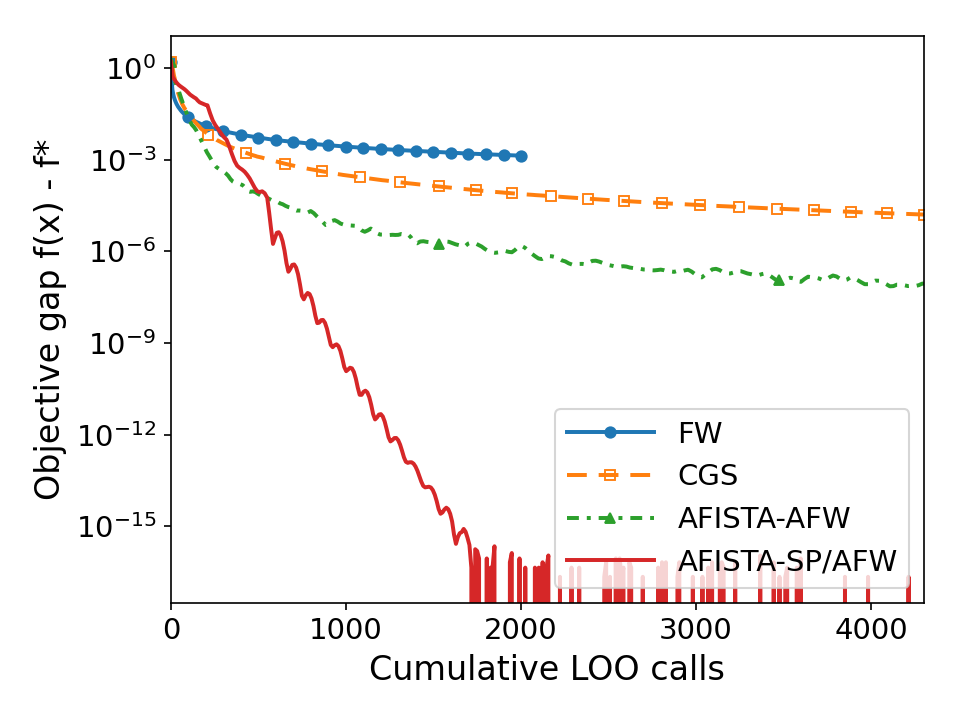}
         \caption*{$r=10, \delta=0.0$}
         \label{fig:y equals x}
     \end{subfigure}
     \hfill
     \begin{subfigure}[b]{0.32\textwidth}
         \centering
         \includegraphics[width=\textwidth]{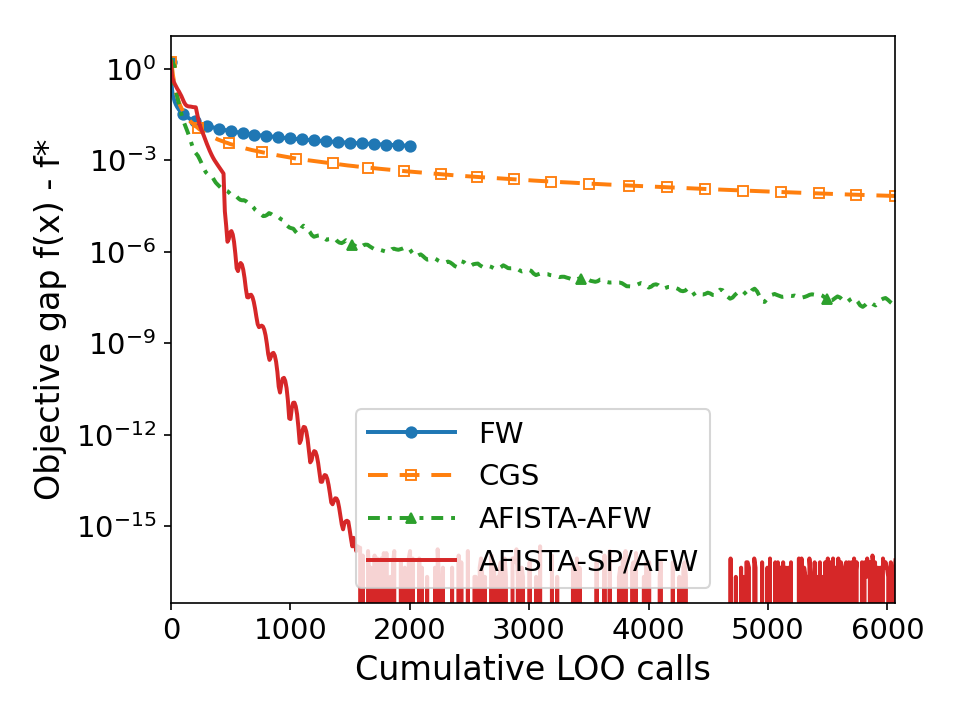}
         \caption*{$r=10, \delta=0.1$}
         \label{fig:three sin x}
     \end{subfigure}
     \hfill
     \begin{subfigure}[b]{0.32\textwidth}
         \centering
         \includegraphics[width=\textwidth]{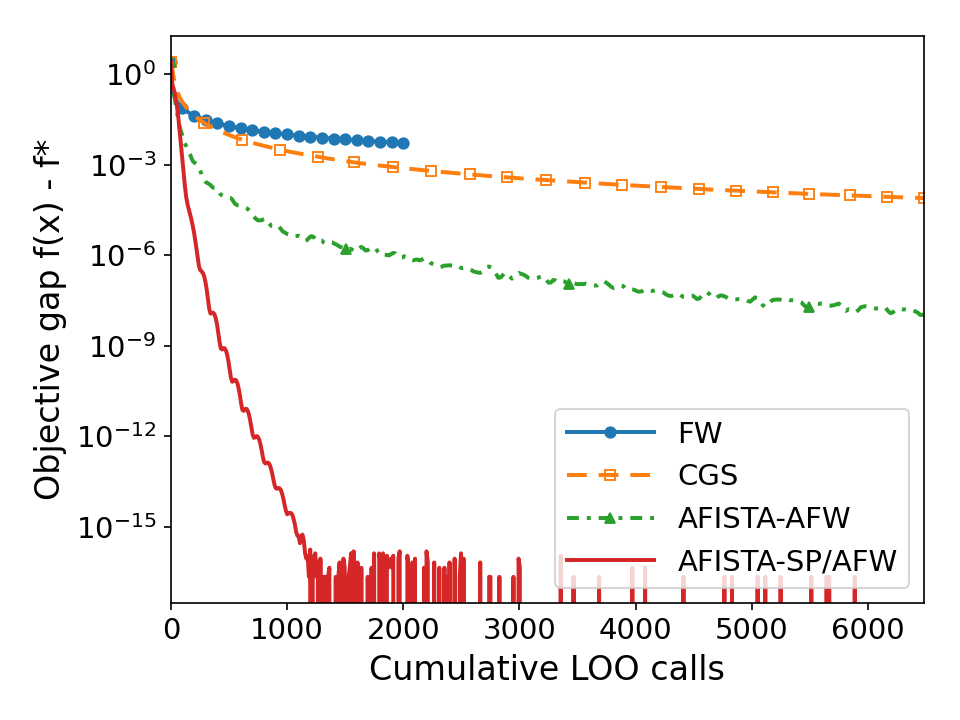}
         \caption*{$r=10, \delta=1.0$}
         \label{fig:five over x}
     \end{subfigure}

     \begin{subfigure}[b]{0.32\textwidth}
         \centering
         \includegraphics[width=\textwidth]{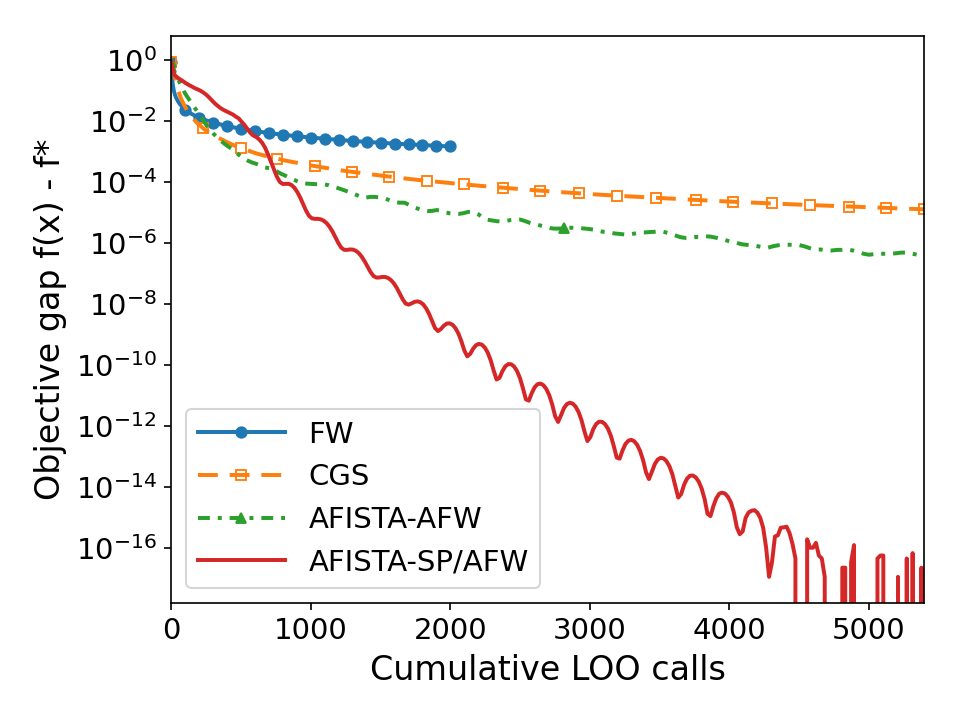}
         \caption*{$r=20, \delta=0.0$}
         \label{fig:y equals x}
     \end{subfigure}
     \hfill
     \begin{subfigure}[b]{0.32\textwidth}
         \centering
         \includegraphics[width=\textwidth]{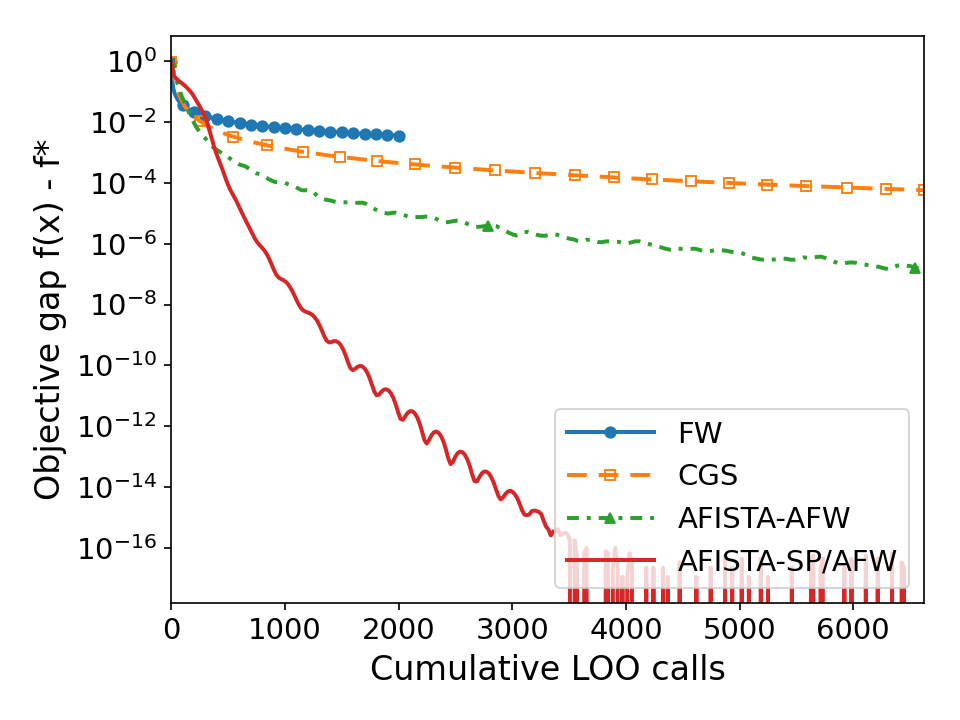}
         \caption*{$r=20, \delta=0.1$}
         \label{fig:three sin x}
     \end{subfigure}
     \hfill
     \begin{subfigure}[b]{0.32\textwidth}
         \centering
         \includegraphics[width=\textwidth]{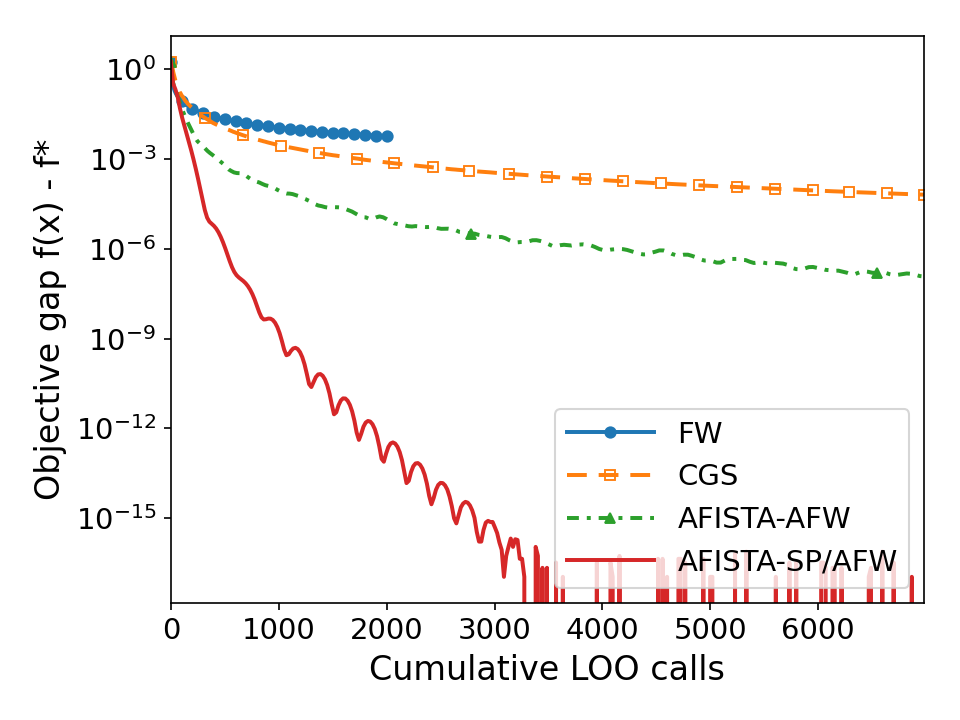}
         \caption*{$r=20, \delta=1.0$}
         \label{fig:five over x}
     \end{subfigure}
     
          \begin{subfigure}[b]{0.32\textwidth}
         \centering
         \includegraphics[width=\textwidth]{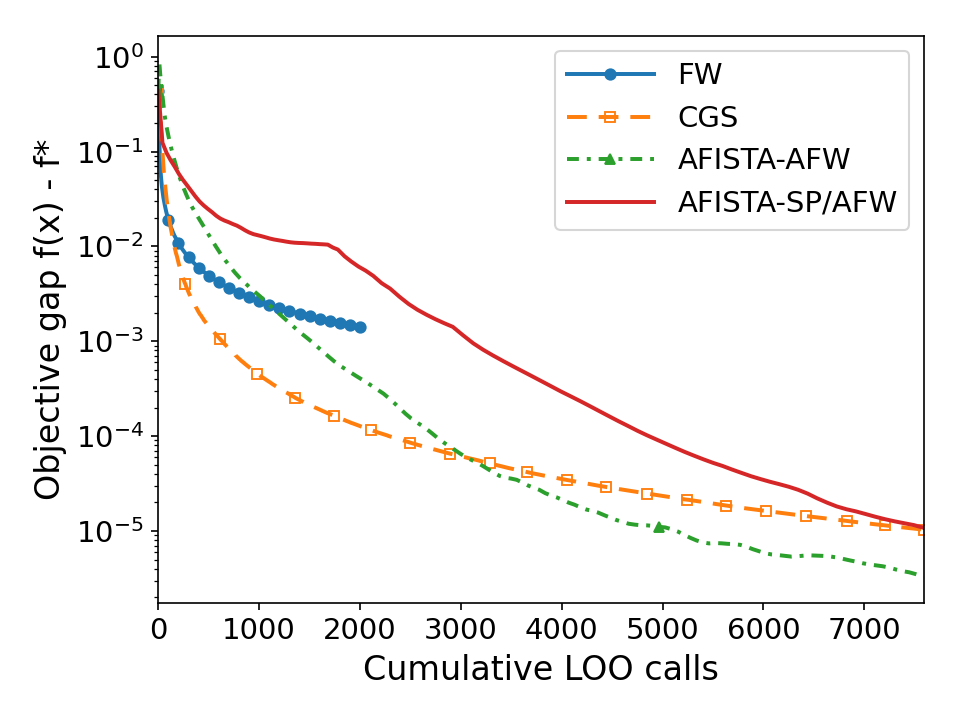}
         \caption*{$r= 40, \delta=0.0$}
         \label{fig:y equals x}
     \end{subfigure}
     \hfill
     \begin{subfigure}[b]{0.32\textwidth}
         \centering
         \includegraphics[width=\textwidth]{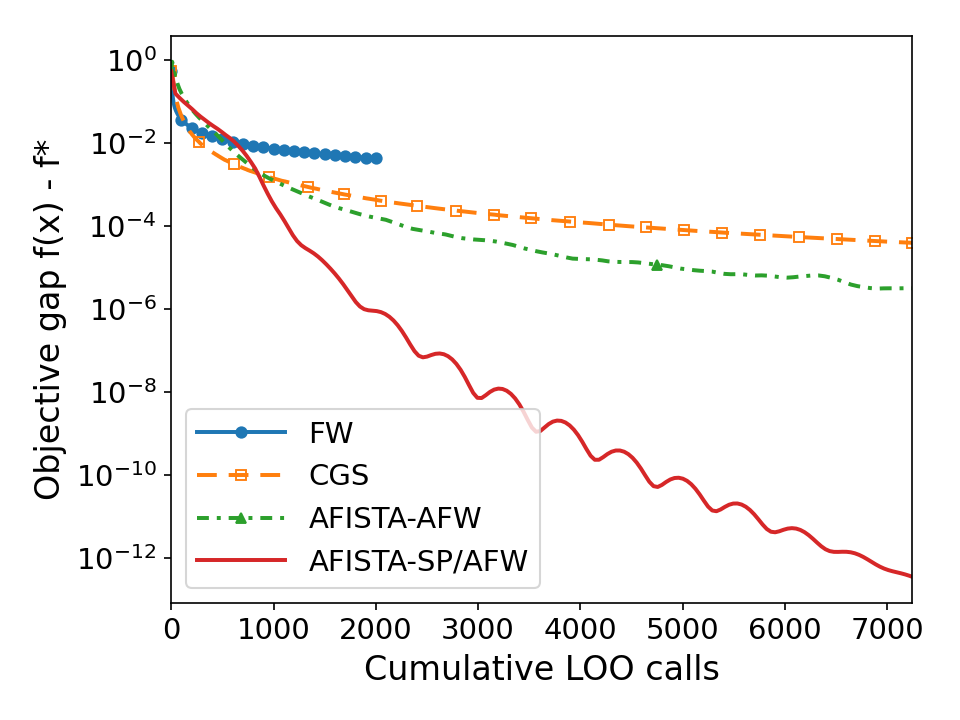}
         \caption*{$r= 40, \delta=0.1$}
         \label{fig:three sin x}
     \end{subfigure}
     \hfill
     \begin{subfigure}[b]{0.32\textwidth}
         \centering
         \includegraphics[width=\textwidth]{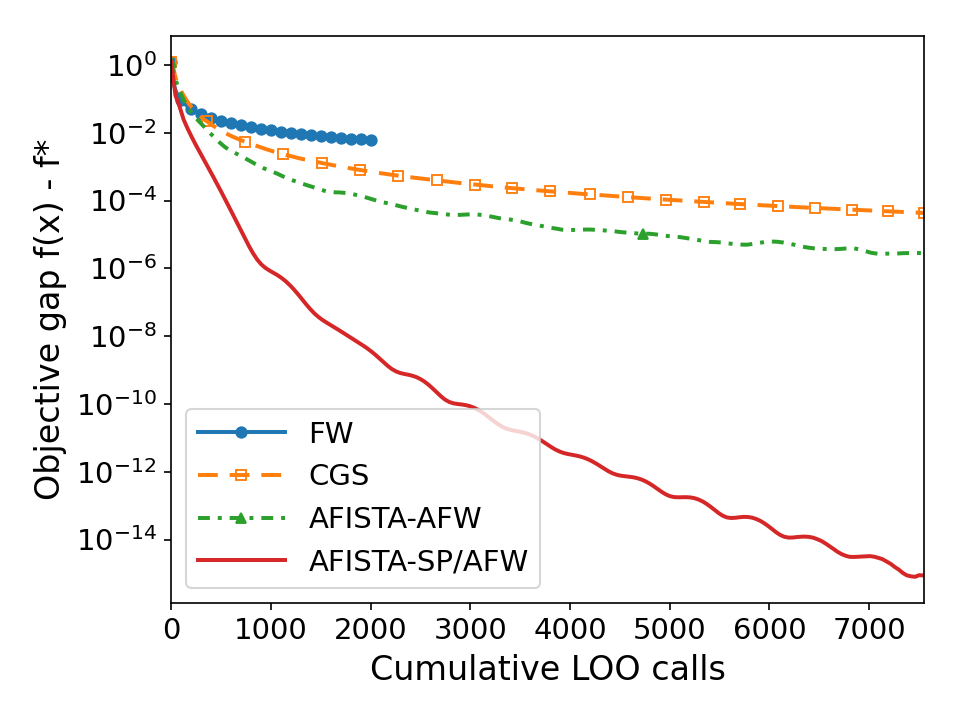}
        \caption*{$r= 40, \delta=1.0$}
         \label{fig:five over x}
     \end{subfigure}

     \begin{subfigure}[b]{0.32\textwidth}
         \centering
         \includegraphics[width=\textwidth]{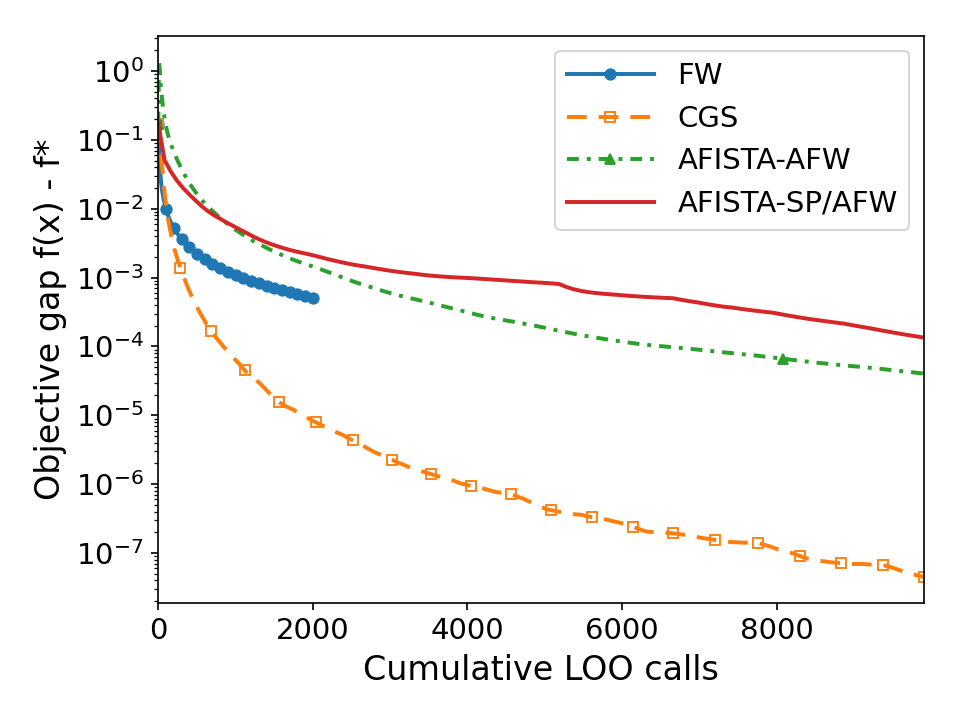}
         \caption*{$r= 80, \delta=0.0$}
         \label{fig:y equals x}
     \end{subfigure}
     \hfill
     \begin{subfigure}[b]{0.32\textwidth}
         \centering
         \includegraphics[width=\textwidth]{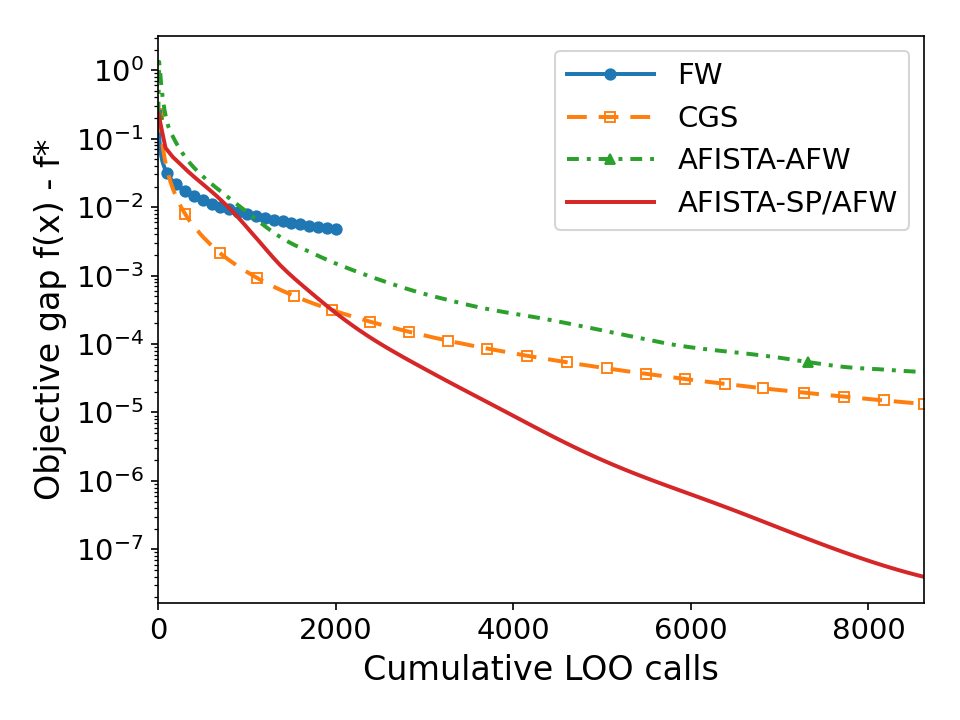}
         \caption*{$r= 80, \delta=0.1$}
         \label{fig:three sin x}
     \end{subfigure}
     \hfill
     \begin{subfigure}[b]{0.32\textwidth}
         \centering
         \includegraphics[width=\textwidth]{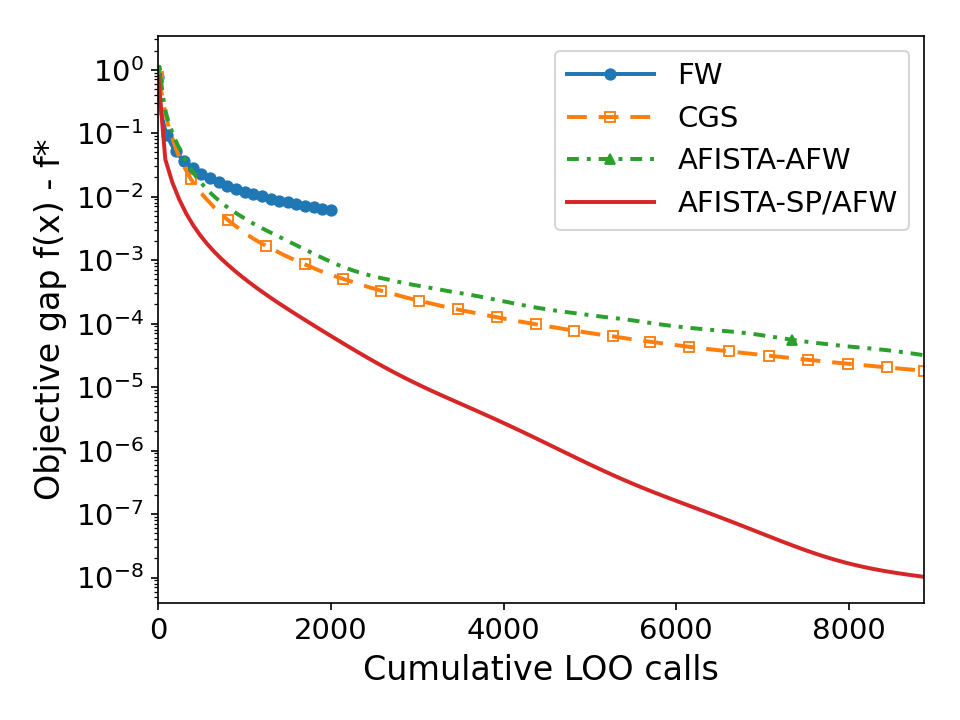}
         \caption*{$r= 80, \delta=1.0$}
     \end{subfigure}
        \caption{Approximation errors vs. number of LOO calls.}
        \label{fig:inner}
\end{figure}

When examining the convergence in terms of LOO calls in Figure \ref{fig:inner}, we see that AFISTA-SP/AFW struggles when $\delta = 0.0$ and $r$ is not very small  compared to the ambient dimension($r\geq 40, n=200$), however this changes dramatically once $\delta >0$. We also see that for small values of $r$, even without strict complementarity our methods can improve significantly over CGS, while CGS has the clear advantage once $r/n$ is large enough ($r=80, n=200$).

\begin{figure}[H]
     \centering
     \begin{subfigure}[b]{0.32\textwidth}
         \centering
         \includegraphics[width=\textwidth]{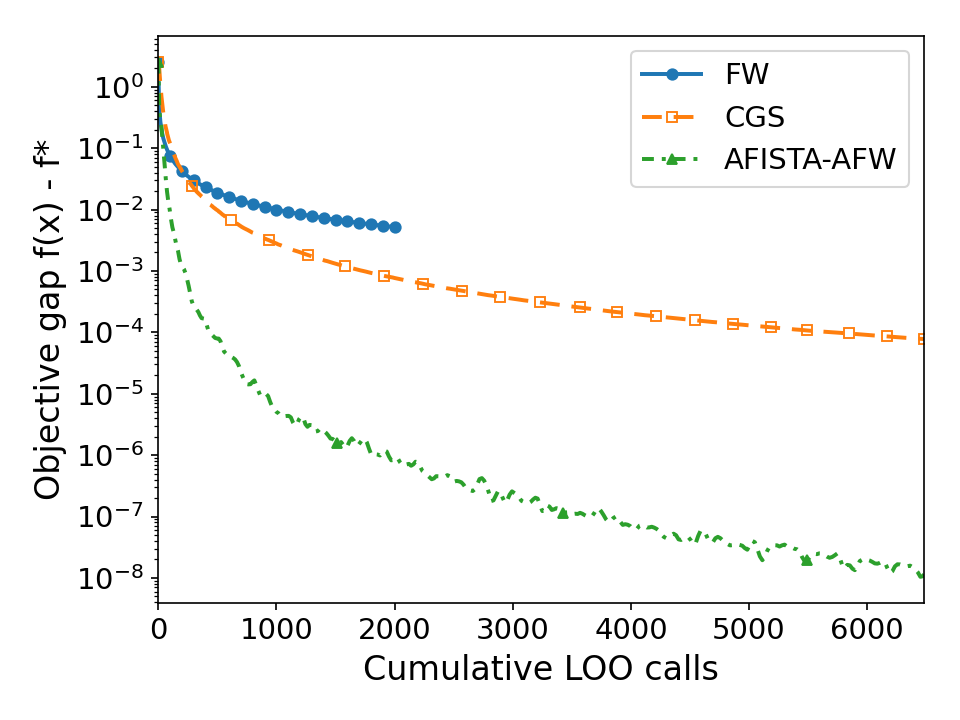}
         \caption*{$r=10, \delta=1.0$}
         \label{fig:y equals x}
     \end{subfigure}
     \hfill
     \begin{subfigure}[b]{0.32\textwidth}
         \centering
         \includegraphics[width=\textwidth]{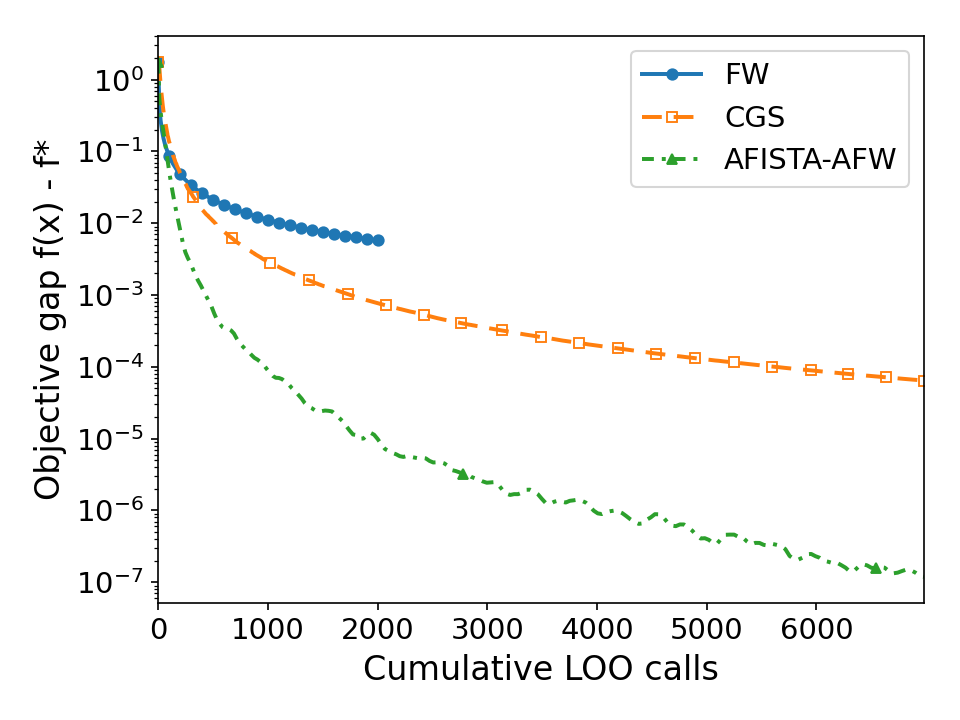}
         \caption*{$r=20, \delta=1.0$}
         \label{fig:three sin x}
     \end{subfigure}
     \hfill
     \begin{subfigure}[b]{0.32\textwidth}
         \centering
         \includegraphics[width=\textwidth]{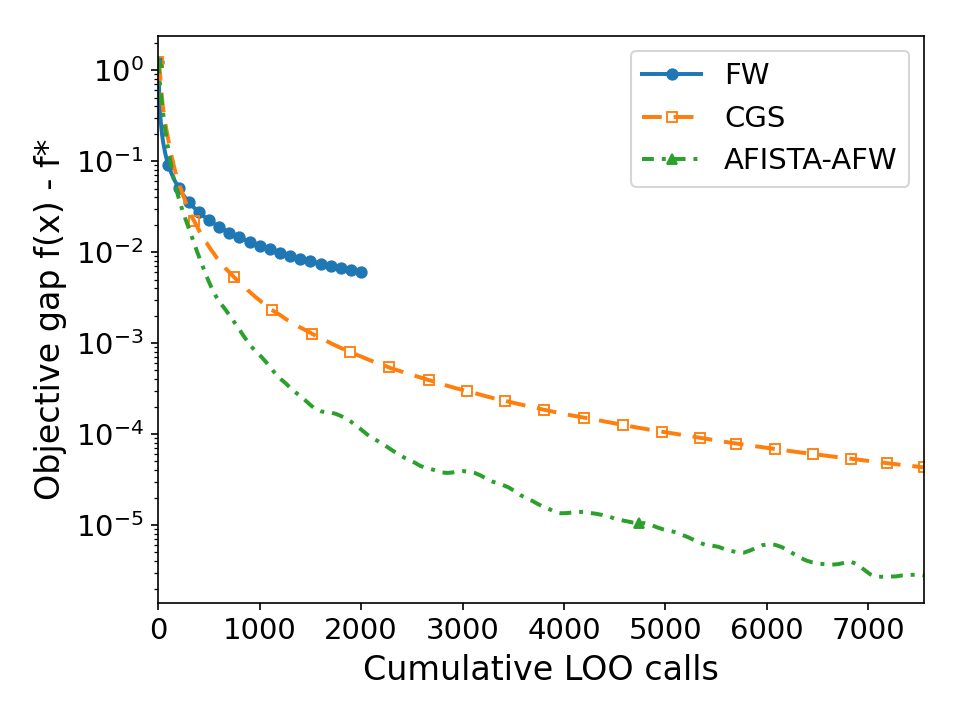}
         \caption*{$r=40, \delta=1.0$}
         \label{fig:five over x}
     \end{subfigure}
        \caption{Approximation errors vs. number of LOO calls without AFISTA-SP/AFW.}
        \label{fig:inner:2}
\end{figure}

When omitting AFISTA-SP/AFW from the comparison (to have a clearer separation of other methods), we see in Figure \ref{fig:inner:2} that AFISTA-AFW indeed benefits significantly from the dimensionality of the optimal face, with a larger margin from CGS for smaller values of $r$.

\section{Acknowledgement}
This work was funded by the European Union (ERC,  ProFreeOpt, 101170791). Views and opinions expressed are however those of the author(s) only and do not necessarily reflect those of the European Union or the European Research Council Executive Agency. Neither the European Union nor the granting authority can be held responsible for them.

\appendix

\section{Approximated FISTA Proofs}
The proofs in this section are closely based on the analysis in \cite{chambolle2015convergence}, with modifications to account for the approximation errors in the sequence $(\x_t)_{t\geq 1}$.\\

The following lemma is a standard argument in the analysis of proximal gradient methods that is adapted to account for approximation errors.
\begin{lemma}\label{lem:FISTAapxProj}
Fix iteration $t\geq 1$ of AFISTA and define the map $\x(\w) := (1-\lambda_t^{-1})\x_{t-1} + \lambda_t^{-1}\w, \w\in\mK$.
For any $\w\in\mK$ we have that 
\begin{align}\label{eq:lem:FISTAapxProj:1}
&\max\left\{f(\x_t^*) + \frac{\beta}{2}\Vert{\x_t^*-\x(\w)}\Vert^2, ~f(\x_t) + \frac{\beta}{2}\Vert{\x_t-\x(\w)}\Vert^2 - \omega_t(\x_t)\right\} - f^* \leq \nonumber \\
&~~~(1-\lambda_t^{-1})f(\x_{t-1}) + \lambda_t^{-1}f(\w) -f^* + \frac{\beta}{2}\Vert{\x(\w)-\y_{t-1}}\Vert^2.
\end{align}
\end{lemma}
\begin{proof}
Fix $\w\in\mK$. From the smoothness of $f$ we have that
\begin{align*}
f(\x_{t}) &\leq f(\y_{t-1}) + \langle{\x_{t} - \y_{t-1}, \nabla{}f(\y_{t-1})}\rangle + \frac{\beta}{2}\Vert{\x_{t}-\y_{t-1}}\Vert^2 \\
& = f(\y_{t-1}) + \phi_t(\x_t) \\
&\underset{(a)}{\leq} f(\y_{t-1}) + \phi_t(\x(\w)) + \langle{\x_t-\x(\w),\nabla\phi_t(\x_t)}\rangle - \frac{\beta}{2}\Vert{\x_t-\x(\w)}\Vert^2 \\
&= f(\y_{t-1}) + \langle{\x(\w) - \y_{t-1}, \nabla{}f(\y_{t-1})}\rangle + \frac{\beta}{2}\Vert{\x(\w)-\y_{t-1}}\Vert^2 \\
&~~~+ \langle{\x_t-\x(\w),\nabla\phi_t(\x_t)}\rangle - \frac{\beta}{2}\Vert{\x_t-\x(\w)}\Vert^2 \\
&\underset{(b)}{\leq} f(\x(\w)) + \frac{\beta}{2}\Vert{\x(\w)-\y_{t-1}}\Vert^2 - \frac{\beta}{2}\Vert{\x_t-\x(\w)}\Vert^2 +\langle{\x_t-\x(\w),\nabla\phi_t(\x_t)}\rangle \\
&\underset{(c)}{\leq} (1-\lambda_t^{-1})f(\x_{t-1}) + \lambda_t^{-1}f(\w) + \frac{\beta}{2}\Vert{\x(\w)-\y_{t-1}}\Vert^2 - \frac{\beta}{2}\Vert{\x_t-\x(\w)}\Vert^2 +\omega_t(\x_t),
\end{align*}
where (a) follows since $\phi_t(\cdot)$ is $\beta$-strongly convex, (b) follows since $f(\cdot)$ is  convex, and (c) follows again from convexity of $f$ and the definition of $\omega_t(\cdot)$.

Subtracting $f^*$ for both sides and rearranging, yields the part of Eq. \eqref{eq:lem:FISTAapxProj:1} which corresponds to the second term inside the max (on the LHS of \eqref{eq:lem:FISTAapxProj:1}).

 To obtain the part of \eqref{eq:lem:FISTAapxProj:1} which corresponds to the first term inside the max, note that if in the above inequalities we replace $\x_t$ with $\x_t^*$, which is the minimizer of $\phi_t$ over $\mK$, due to the first-order optimality condition, the term $\langle{\x_t^*-\x(\w),\nabla\phi_t(\x_t^*)}\rangle$ is not positive and can be omitted.

\end{proof}

\begin{lemma}\label{lem:FISTA:funcConv}
Consider Algorithm AFISTA with $a \geq 2$ and 
suppose that for all $t$, $\omega_t(\x_t) \leq \nu_t$ for some non-negative sequence $(\nu_t)_{t\geq 1}$. Then, for all $T\geq 0$ and any $\x^*\in\mX^*$ it holds that,
\begin{align}\label{lem:FISTA:funcConv:Res1}
h_{T+1} \leq \frac{1}{\lambda_{T+1}^2}\left({\frac{\beta}{2}\Vert{\x_0-\x^*}\Vert^2 + \sum_{t=0}^T\lambda_{t+1}^2\nu_{t+1}}\right).
\end{align}
Furthermore, denoting $\rho_t = \lambda_{t-1}^2 - \lambda_{t}^2+\lambda_{t}$ for all $t\geq 2$, it holds that
\begin{align}\label{lem:FISTA:funcConv:Res2}
\sum_{t=1}^T\rho_{t+1}h_t \leq \frac{\beta}{2}\Vert{\x_0-\x^*}\Vert^2 + \sum_{t=0}^{T-1}\lambda_{t+1}^2\nu_{t+1}.
\end{align}
\end{lemma}
\begin{proof}
Denote  $\u_t = \x_{t-1} + \lambda_t(\x_{t} - \x_{t-1})$ for all $t\geq 1$ and $\u_0 = \x_0$.

Fix some $t\geq 0$. Applying Lemma \ref{lem:FISTAapxProj} for iteration $t+1$ and $\w=\x^*$, and denoting $\x = (1-\lambda_{t+1}^{-1})\x_t + \lambda_{t+1}^{-1}\x^*$, we have that
\begin{align}\label{eq:lem:FISTAvalConv:1}
&h_{t+1} + \frac{\beta}{2}\Vert{\x_{t+1}-\x}\Vert^2 \leq (1-\lambda_{t+1}^{-1})h_t  + \frac{\beta}{2}\Vert{(1-\lambda_{t+1}^{-1})\x_t + \lambda_{t+1}^{-1}\x^*-\y_{t}}\Vert^2 + \nu_{t+1}.
\end{align}
For $t\geq 1$, using the definition of $\y_t$ in Eq. \eqref{eq:FISTA:yup}, we have that
\begin{align*}
h_{t+1} + \frac{\beta}{2}\Vert{\x_{t+1}-\x}\Vert^2 &\leq  (1-\lambda_{t+1}^{-1})h_t + \frac{\beta}{2}\Vert{\lambda_{t+1}^{-1}\left({\x^*-\x_{t-1}-\lambda_t(\x_t-\x_{t-1})}\right)}\Vert^2 + \nu_{t+1}\nonumber\\
&=  (1-\lambda_{t+1}^{-1})h_t + \frac{\beta}{2\lambda_{t+1}^2}\Vert{\x^*-\u_{t}}\Vert^2 + \nu_{t+1},
\end{align*}
and observing that
\begin{align*}
\x_{t+1}-\x = \lambda_{t+1}^{-1}\left({\x_t + \lambda_{t+1}(\x_{t+1}-\x_t) - \x^*}\right) = \lambda_{t+1}^{-1}(\u_{t+1} - \x^*),
\end{align*}
we have that for all $t\geq 1$,
\begin{align*}
&h_{t+1} - \left({1-\lambda_{t+1}^{-1}}\right)h_t \leq \frac{\beta}{2\lambda_{t+1}^2}\Vert{\u_t-\x^*}\Vert^2 - \frac{\beta}{2\lambda_{t+1}^2}\Vert{\u_{t+1}-\x^*}\Vert_2^2 + \nu_{t+1}.
\end{align*}
Also, for $t=0$, starting from Eq. \eqref{eq:lem:FISTAvalConv:1} and recalling that $\y_0 = \x_0=\u_0$, $\lambda_1 = 1$, and so $\x_1 = \u_1$ and $\x = \x^*$, we have that
\begin{align*}
h_1 + \frac{\beta}{2\lambda_1^2}\Vert{\u_1 - \x^*}\Vert^2 \leq (1-\lambda_1^{-1})h_0 + \frac{\beta}{2\lambda_1^2}\Vert{\u_0-\x^*}\Vert^2 +\nu_1.
\end{align*}

Thus, by rearranging the last two inequalities we have that for all $t\geq 0$,
\begin{align}\label{eq:lem:FISTAvalConv:2}
\lambda_{t+1}^2h_{t+1} - (\lambda_{t+1}^2 - \lambda_{t+1})h_t &\leq \frac{\beta}{2}\left({\Vert{\u_t-\x^*}\Vert^2 - \Vert{\u_{t+1}-\x^*}\Vert^2}\right)  + \lambda_{t+1}^2\nu_{t+1}.
\end{align}
Summing the above from $t=0$ to $T$,  using the definition of $\rho_t$ in the lemma and recalling that $\u_0 = \x_0$ and $\lambda_1=1$, gives that for all $T\geq 0$,
\begin{align*}
\lambda_{T+1}^2h_{T+1} + \sum_{t=1}^T\rho_{t+1}h_t &\leq \frac{\beta}{2}\Vert{\x_0-\x^*}\Vert^2 + \sum_{t=0}^T\lambda_{t+1}^2\nu_{t+1}.
\end{align*}

The first part of the lemma (Eq. \eqref{lem:FISTA:funcConv:Res1}) follows from the choice of the sequence $(\lambda_t)_{t\geq 1}$  and the observation that it implies that $\rho_t \geq 0$. Indeed, a simple calculation yields:
\begin{align}\label{eq:thm:funValConv:1}
\rho_t &= \frac{1}{a^2}\left({(t+a-2)^2 - (t+a-1)^2 + a(t+a-1)}\right) \nonumber \\
&= \frac{1}{a^2}\left({(a-2)t+a^2-3a+3}\right) \geq 0,
\end{align}
where the last inequality is due to the assumption $a \geq 2$.

We continue to prove the second part of the lemma.
Applying Lemma \ref{lem:FISTAapxProj} again with $t=T+1$ and $\w = \x^*$, and denoting $\x = (1-\lambda_{T+1}^{-1})\x_T + \lambda_{T+1}^{-1}\x^*$, $h_{T+1}^* = f(\x_{T+1}^*) - f(\x^*)$, we have that
\begin{align*}
h_{T+1}^* + \frac{\beta}{2}\Vert{\x_{T+1}^*-\x}\Vert^2 &\leq (1-\lambda_{T+1}^{-1})h_T  + \frac{\beta}{2}\Vert{(1-\lambda_{T+1}^{-1})\x_T + \lambda_{T+1}^{-1}\x^*-\y_{T}}\Vert^2 \\
&=  (1-\lambda_{T+1}^{-1})h_T + \frac{\beta}{2}\Vert{\lambda_{T+1}^{-1}\left({\x^*-\x_{T}-(\lambda_T-1)(\x_T-\x_{T-1})}\right)}\Vert^2\nonumber\\
&=  (1-\lambda_{T+1}^{-1})h_T + \frac{\beta}{2\lambda_{T+1}^2}\Vert{\x^*-\u_{T}}\Vert^2,
\end{align*}
which by denoting $\u_{T+1}^* = \x_T + \lambda_{T+1}(\x_{T+1}^*-\x_T)$ gives,
\begin{align*}
\lambda_{T+1}^2h_{T+1}^* - (\lambda_{T+1}^2 - \lambda_{T+1})h_T &\leq \frac{\beta}{2}\left({\Vert{\u_T-\x^*}\Vert^2 - \Vert{\u_{T+1}^*-\x^*}\Vert^2}\right).
\end{align*}
Summing Eq. \eqref{eq:lem:FISTAvalConv:2}  from $t=0$ to $T-1$ and adding the above inequality gives,
\begin{align*}
\lambda_{T+1}^2h_{T+1}^* + \sum_{t=1}^T\rho_{t+1}h_t &\leq \frac{\beta}{2}\Vert{\x_0-\x^*}\Vert^2 + \sum_{t=0}^{T-1}\lambda_{t+1}^2\nu_{t+1},
\end{align*}
which, due to the non-negativity of $h_{T+1}^*$, yields the second part of the lemma (Eq. \eqref{lem:FISTA:funcConv:Res2}).
\end{proof}

\begin{lemma}\label{lem:FISTA:distConv}
Consider Algorithm AFISTA with $a \geq 5$ and suppose that for all $t$, $\omega_t(\x_t) \leq \nu_t$ for some non-negative sequence $(\nu_t)_{t\geq 1}$. Then, for all $T\geq 0$ and any $\x^*\in\mX^*$ it holds that,
\begin{align*}
d_{T+1}^* &\leq  \frac{1}{\lambda_{T+1}^2}\left({\Vert{\x_0-\x^*}\Vert^2 + \frac{3}{\beta}\sum_{t=0}^{T-1}\lambda_{t+1}^2\nu_{t+1}}\right),\\
d_{T+1} &\leq  \frac{1}{\lambda_{T+1}^2}\left({\Vert{\x_0-\x^*}\Vert^2 + \frac{3}{\beta}\sum_{t=0}^{T}\lambda_{t+1}^2\nu_{t+1}}\right).
\end{align*}
\end{lemma}
\begin{proof}
We first prove the bound on $d_{T+1}^*$ and then on $d_{T+1}$.
Fix  $t\geq 0$. Applying Lemma \ref{lem:FISTAapxProj} for iteration $t+1$ and with $\w=\x_t$ gives,
\begin{align*}
h_{t+1} + \beta{}d_{t+1} &\leq h_t + \frac{\beta}{2}\Vert{\x_t-\y_t}\Vert^2 + \nu_{t+1}. 
\end{align*}
For $t\geq 1$, using the definition of $\y_t$ in Eq. \eqref{eq:FISTA:yup}, we have that
\begin{align*}
h_{t+1} + \beta{}d_{t+1} \leq h_t + \beta\left({\frac{\lambda_t-1}{\lambda_{t+1}}}\right)^2d_t + \nu_{t+1}.
\end{align*}
Similarly, denoting $h_{T+1}^* = f(\x_{T+1}^*) - f(\x^*)$, and using Lemma \ref{lem:FISTAapxProj} for iteration $T+1$ and with $\w=\x_T$, we also have that
\begin{align*}
h_{T+1}^* + \beta{}d_{T+1}^* \leq h_T +  \beta\left({\frac{\lambda_T-1}{\lambda_{T+1}}}\right)^2d_T.
\end{align*}
Denoting $\theta_t = \frac{\lambda_t -1}{\lambda_{t+1}}$, these yield
\begin{align}
d_{t+1} - \theta_t^2d_t &\leq \frac{1}{\beta}\left({h_t - h_{t+1}}\right) + \frac{\nu_{t+1}}{\beta} \quad\forall t\geq 1; \label{eq:FISTA:distConv:1}\\
d_{T+1}^* - \theta_{T}^2d_T &\leq \frac{1}{\beta}\left({h_T - h_{T+1}^*}\right). \label{eq:FISTA:distConv:2}
\end{align}

Multiplying \eqref{eq:FISTA:distConv:1} by $(t+a)^2$ on both sides and summing from $t=1$ to $t=T-1$, and then adding to it \eqref{eq:FISTA:distConv:2} multiplied by $(T+a)^2$ on both sides, we obtain
\begin{align*}
&(T+a)^2\left({d_{T+1}^* - \theta_T^2d_T}\right) +  \sum_{t=1}^{T-1}(t+a)^2\left({d_{t+1} - \theta_t^2d_t}\right) \\
&\leq  \frac{1}{\beta}(T+a)^2\left({h_T - h_{T+1}^*}\right)+ \frac{1}{\beta}\sum_{t=1}^{T-1}(t+a)^2\left({h_t - h_{t+1}}\right) + \frac{1}{\beta}\sum_{t=1}^{T-1}(t+a)^2\nu_{t+1},
\end{align*}
which simplifies to (recall $\theta_1 = 0$),
\begin{align*}
&(T+a)^2d_{T+1}^* + \sum_{t=2}^T\left({(t+a-1)^2-(t+a)^2\theta_t^2}\right)d_t \\
& \leq \frac{1}{\beta}\left({(1+a)^2h_1 + \sum_{t=2}^{T}\left({(t+a)^2-(t+a-1)^2}\right)h_t}\right) + \frac{1}{\beta}\sum_{t=1}^{T-1}(t+a)^2\nu_{t+1}.
\end{align*}
Using $\theta_t = \frac{\lambda_{t}-1}{\lambda_{t+1}} = \frac{t-1}{t+a}$, the above  further simplifies to
\begin{align}\label{eq:thm:itDist:1}
(T+a)^2d_{T+1}^* \leq \frac{1}{\beta}\left({(1+a)^2h_1 + \sum_{t=2}^T(2t+2a-1)h_t}\right) + \frac{1}{\beta}\sum_{t=1}^{T-1}(t+a)^2\nu_{t+1}.
\end{align}

Lemma \ref{lem:FISTA:funcConv} implies that
\begin{align*}
\sum_{t=1}^T\rho_{t+1}h_t \leq \frac{\beta}{2}\Vert{\x_0-\x^*}\Vert^2 + \sum_{t=0}^{T-1}\lambda_{t+1}^2\nu_{t+1},
\end{align*}
which by Eq. \eqref{eq:thm:funValConv:1} and the choice of the sequence $(\lambda_t)_{t\geq 1}$ further implies that,
\begin{align*}
\sum_{t=1}^T\frac{1}{a^2}\left({(a-2)t+a^2-2a+1}\right)h_t &\leq \frac{\beta}{2}\Vert{\x_0-\x^*}\Vert^2 + \sum_{t=0}^{T-1}\frac{(t+a)^2}{a^2}\nu_{t+1}.
\end{align*}
In particular, a simple calculation verifies that for $a \geq 5$ we have that,
\begin{align*}
\frac{1}{2a^2}\left({(1+a)^2h_1 + \sum_{t=2}^T\left({2t+2a-1}\right)h_t}\right) &\leq \frac{\beta}{2}\Vert{\x_0-\x^*}\Vert^2 + \sum_{t=0}^{T-1}\frac{(t+a)^2}{a^2}\nu_{t+1}.
\end{align*}
Plugging this inequality into \eqref{eq:thm:itDist:1} and simplifying we obtain,
\begin{align*}
d_{T+1}^* \leq \frac{1}{(T+a)^2}\left({a^2\Vert{\x_0-\x^*}\Vert^2 + \frac{3}{\beta}\sum_{t=0}^{T-1}(t+a)^2\nu_{t+1}}\right).
\end{align*}
Plugging-in the definition of $\lambda_t$ for all $t$ into the above inequality and rearranging, yields the bound on $d_{T+1}^*$.

To prove the second part of the lemma, the upper-bound on $d_{T+1}$, we go back to Eq. \eqref{eq:FISTA:distConv:1} and \eqref{eq:FISTA:distConv:2}, but this time, we shall only use \eqref{eq:FISTA:distConv:1}, i.e., we shall multiply it on both sides by $(t+a)^2$ and sum from $t=1$ to $t=T$. This will yield the inequality
\begin{align}\label{eq:thm:itDist:1b}
(T+a)^2d_{T+1} \leq \frac{1}{\beta}\left({(1+a)^2h_1 + \sum_{t=2}^T(2t+2a-1)h_t}\right) + \frac{1}{\beta}\sum_{t=1}^{T}(t+a)^2\nu_{t+1},
\end{align}
instead of the previous Eq. \eqref{eq:thm:itDist:1} (note that now the sum on $\nu_{t+1}$ in the RHS include an additional term --- $\nu_{T+1}$). From here we continue exactly as in the derivation following Eq. \eqref{eq:thm:itDist:1b} above and we shall obtain the bound 
\begin{align*}
d_{T+1} \leq \frac{1}{(T+a)^2}\left({a^2\Vert{\x_0-\x^*}\Vert^2 + \frac{3}{\beta}\sum_{t=0}^{T}(t+a)^2\nu_{t+1}}\right),
\end{align*}
and the result follows again from plugging-in the definition of $\lambda_t$ for all $t$ into the above inequality and rearranging 
\end{proof}

\section{Lemma 7}

\begin{lemma}\label{lem:smoothIneq}
Let $\Psi:\mathbb{E}\to\mathbb{R}$ be convex and $\beta_{\Psi}$-smooth,
and let $\mK\subseteq\mathbb{E}$ be convex and compact.
The gradient $\nabla\Psi$ is constant over
$\arg\min_{\z\in\mK}\Psi(\z)$, and, for any $\z\in\mK$ and any
$\z^*\in\arg\min_{\z\in\mK}\Psi(\z)$,
\begin{equation*}
    \|\nabla\Psi(\z)-\nabla\Psi(\z^*)\|^2
    \leq
    2\beta_{\Psi}\bigl(\Psi(\z)-\Psi(\z^*)\bigr).
\end{equation*}
\end{lemma}

\begin{proof}
A standard inequality for smooth convex functions says that  for any $\x,\y\in\mathbb{E}$,
\begin{equation}\label{eq:smooth-convex-gradient}
    \Psi(\x)
    \geq
    \Psi(\y)
    + \langle{\x-\y, \nabla\Psi(\y)}\rangle
    + \frac{1}{2\beta_{\Psi}}
      \|\nabla\Psi(\x)-\nabla\Psi(\y)\|^2 .
\end{equation}

Let $\z_1^*,\z_2^*\in\arg\min_{\z\in\mK}\Psi(\z)$. By the first-order optimality condition,
$\langle{\z_2^*-\z_1^*, \nabla\Psi(\z_1^*)}\rangle \geq 0$.
Applying~\eqref{eq:smooth-convex-gradient} with
$\x=\z_2^*$ and $\y=\z_1^*$ then gives,
\[
    0
    \geq
    \frac{1}{2\beta_{\Psi}}
      \|\nabla\Psi(\z_2^*)-\nabla\Psi(\z_1^*)\|^2,
\]
and so $\nabla\Psi(\z_1^*)=\nabla\Psi(\z_2^*)$, which proves the first part.

Now fix $\z\in\mK$ and
$\z^*\in\arg\min_{\z\in\mK}\Psi(\z)$. Again, using the first-order
optimality condition $\langle{\z-\z^*,\nabla\Psi(\z^*)}\rangle\geq0$, and applying~\eqref{eq:smooth-convex-gradient} with
$\x=\z$ and $\y=\z^*$  yields
\[
    \Psi(\z)-\Psi(\z^*)
    \geq
    \frac{1}{2\beta_{\Psi}}
    \|\nabla\Psi(\z)-\nabla\Psi(\z^*)\|^2,
\]
which proves the second part.
\end{proof}

\section{Proof of Theorem \ref{thm:afw}}

Before proving the theorem we need the following lemma which is adapted from Lemma 5.5 in \cite{garber2016linearly}. A similar adaptation is also given in \cite{garber2020revisiting}, however for completeness and clarity of presentation we detail it here.
\begin{lemma}\label{lem:poly:dist}
Let $\mP\subset\reals^n$ be a polytope of the form $\{\x\in\reals^n~|~\tilde{\A}_1\x = \tilde{\b}_1, \tilde{\A}_2\x \leq \tilde{\b}_2\}$, $\tilde{\A}_1\in\reals^{\tilde{m}_1\times n}$, $\tilde{\A}_2\in\reals^{\tilde{m}_2\times n}$, such that $\mathrm{aff}(\mP) = \{\x ~|~\tilde{\A}_1\x = \tilde{\b}_1\}$, and let $\mV_{\mP}$ denote its set of vertices. Let $\dim\mP$ be as defined in \eqref{eq:polyDim} and let $\mu_{\mP}$ be as defined in \eqref{eq:mu}. Fix some $\x\in\mP$ given as a convex combination $\x = \sum_{j=1}^m\rho_j\v_j$, where $\rho_j > 0$ for all $j$ and $\sum_{j=1}^m\rho_j =1$, and $\{\v_1,\dots,\v_m\}\subseteq\mV_{\mP}$. Then, for any $\y\in\mP$ there exists some $\z\in\mP$ and scalars $\Delta_1,\dots,\Delta_m$ satisfying $\Delta_j \in [0, \rho_j]$  for all $j$, $\sum_{j=1}^m\Delta_j \leq \mu_{\mP}\sqrt{\dim\mP}\Vert{\x-\y}\Vert$, such that $\y$ can be written as $\y = \sum_{j=1}^m(\rho_j-\Delta_j)\v_j + \sum_{j=1}^m\Delta_j\z$.
\end{lemma}
\begin{proof}
Since the lemma and its proof are a simple refinement of Lemma 5.5 in \cite{garber2016linearly} (a similar refinement was also used in \cite{garber2020revisiting}, Lemma 2), we only detail the simple differences. Lemma 5.5 in  \cite{garber2016linearly} shows the result of our lemma holds but with the ambient dimension $n$ instead of $\dim\mP$, that is Lemma 5.5 in \cite{garber2016linearly} establishes that:
\begin{align}\label{eq:lem:poly:dist:1}
\sum_{j=1}^m\Delta_j \leq \mu_{\mP}\sqrt{n}\Vert{\x-\y}\Vert.
\end{align}

In the sequel, for any matrix $\A\in\reals^{\tilde{m}\times n}$ and $i\in[\tilde{m}]$ we shall denote by $\A(i)$ the column vector corresponding to the $i$th row of $\A$.

The dependence on $n$ in the RHS of \eqref{eq:lem:poly:dist:1} comes from an upper-bound on the cardinality of a set $C_0(\z) \subseteq [\tilde{m}_2]$, where $C_0(\z)$ is any subset of $[\tilde{m}_2]$ (i.e., $C_0$ indexes inequality constraints defining the polytope $\mP$) satisfying the following conditions:
\begin{enumerate}
\item
the vectors $\{\tilde{\A}_2(i)\}_{i\in{}C_0(\z)}\}$ are linearly independent;
\item
$\langle{\tilde{\A}_2(i),\z}\rangle = \tilde{\b}_2(i)$ for all $i\in{}C_0(\z)$;
\item
for any $j\in[m]$ with $\Delta_j > 0$ there exists $i_j\in{}C_0(\z)$ such that $\langle{\tilde{\A}_2(i_j), \v_j}\rangle < \tilde{\b}_2(i_j)$.
\end{enumerate}
Indeed the bound $\vert{C_0(\z)}\vert \leq n$ holds trivially due to the first condition. 

We now show however a refined bound that only scales with $\dim\mP$. Let $C_0(\z)\subseteq[\tilde{m}_2]$ be a set of minimal cardinality satisfying the above three conditions. First note that this implies that for any $i\in{}C_0(\z)$ there must exist some $j_i\in[m]$ such that $\v_{j_i}$ (a vertex in the convex sum yielding $\x$, as assumed in the lemma) satisfies with equality all inequality constraints indexed by $C_0(\z)$ except for constraint number $i$. We shall say that constraint $i$ is \textit{critical} for $\v_{j_i}$. If this is not the case for some $i\in{}C_0(\z)$, then it is redundant  in $C_0(\z)$ (i.e., removing it won't violate any of the three conditions above), which contradicts the minimal cardinality of $C_0(\z)$. 

We now argue that it must hold that for any $i\in{}C_0(\z)$, the vector $\tilde{\A}_2(i)$ is linearly independent of $\{\tilde{\A}_2(k)\}_{k\in{}C_0(\z)\setminus\{i\}}\}\cup\{\tilde{\A}_1(k)\}_{k\in[\tilde{m}_1]}$. To see why this is true, suppose by way of contradiction that this does not hold for some $i\in{}C_0(\z)$ and let $\v_{j_i}$ be a vertex for which $i$ is a critical constraint. Observe that $\v_{j_i}$ must satisfy all equality constraints $\tilde{\A}_1\v_{j_i}=\tilde{\b}_1$ and also the inequality constraints indexed by $C_0(\z)\setminus\{i\}$. Thus, if $\tilde{\A}_2(i)$ is linearly dependent on $\{\tilde{\A}_2(k)\}_{k\in{}C_0(\z)\setminus\{i\}}\}\cup\{\tilde{\A}_1(k)\}_{k\in[\tilde{m}_1]}$, it must follow that $\v_{j_i}$ also satisfies with equality the inequality constraint $i$, which leads to a contradiction.

Thus, since  for any $i\in{}C_0(\z)$ the vector $\tilde{\A}_2(i)$ is linearly independent of $\{\tilde{\A}_2(k)\}_{k\in{}C_0(\z)\setminus\{i\}}\}\cup\{\tilde{\A}_1(k)\}_{k\in[\tilde{m}_1]}$, we indeed have that 
\begin{align*}
\vert{C_0(\z)}\vert \leq n - \dim\textrm{span}\left({\{\tilde{\A}_1(i)\}_{i\in[\tilde{m}_1]}}\right) = \dim\mP.
\end{align*}
\end{proof}

We can now prove Theorem  \ref{thm:afw}.
\begin{proof}[Proof of Theorem  \ref{thm:afw}]
First, we note that the claim that $\omega_t(\x_t) \leq \nu_t$ follows from an immediate application of Lemma \ref{lem:slide} and the stopping condition of the algorithm.
 
For all $i\in\mathbb{N}_+$ denote $g_i = \Phi_t(\w_i) - \Phi_t(\w^*)$, where $\w^*=\arg\min_{\w\in\mK}\Phi_t(\w)$. Additionally, denote the dual gap on iteration $i$, $q_i = \max_{\u\in\mK}\langle{\w_i - \u, \nabla\Phi_t(\w_i)}\rangle$. Note that due to convexity of $\Phi_t(\cdot)$, we have that $g_i \leq q_i$ for all $i$.

Let $\beta_t = \beta/\lambda_t^2$ denote the smoothness parameter of $\Phi_t(\cdot)$.
On each iteration $i$ which is not a drop step, i.e., not a step in which the away direction is chosen and $\gamma_i = \gamma_{\max}$, we have using the smoothness of $\Phi_t(\cdot)$ that 
\begin{align}\label{eq:afw:1}
\forall \gamma\in[0,1]: \quad \Phi_t(\w_{i+1}) \leq \Phi_t(\w_i) - \gamma{}q_i + \frac{\gamma^2\beta_tD^2}{2},
\end{align}
see for instance the very short proof of Lemma 1 in \cite{garber2020revisiting}. 
This is the exact single iteration error reduction as in the standard Frank-Wolfe algorithm with line-search \cite{Jaggi13}. Thus, the same convergence argument for the sequence of dual gaps $(q_i)_{i\geq 1}$ (Theorem 2 in   \cite{Jaggi13}) holds here with a single distinction: for Algorithm \ref{alg:awayfw} we only count the iterations that are not drop steps (importantly drop steps cannot increase the objective $\Phi_t$). Since for any number of iterations $\tau$, the number of drop steps after $\tau$ iterations cannot exceed $(\tau+1)/2$ (see Observation 1 in  \cite{garber2020revisiting}), we have that the dual convergence rate $\min_{1\leq i \leq \tau}q_i = O(\beta_tD^2/\tau)$ of the standard Frank-Wolfe  algorithm (Theorem 2 in   \cite{Jaggi13}), holds also for Algorithm \ref{alg:awayfw}. Plugging-in the value of $\beta_t$ and the stopping condition of the algorithm ($q_i \leq \nu_t$), this proves the first term inside the $\min$ in Eq. \eqref{eq:afw:conv}.\\

For the second term inside the $\min$ in Eq. \eqref{eq:afw:conv}, we adapt the linear convergence argument from \cite{garber2020revisiting} (Theorem 5) which is as follows.
Consider  some iteration $i$ and write $\w_i$ as a convex combination of vertices in $\mV$, i.e.,
$\w_i = \sum_{j=1}^m\rho_j\v_j$, $\{\v_j\}_{j\in[m]}\subseteq\mV$, $\rho_j> 0~\forall j$, $\sum_{j=1}^m\rho_j=1$. Suppose without loss of generality that $\v_1,\dots,\v_m$ are ordered such that  $\langle{\v_1,\nabla{}\Phi_t(\w_i)}\rangle \geq  \langle{\v_2,\nabla{}\Phi_t(\w_i)}\rangle\geq \dots \geq \langle{\v_m,\nabla{}\Phi_t(\w_i)}\rangle$. 

According to Lemma \ref{lem:poly:dist}, there exist scalars $\Delta_1,\dots,\Delta_m$ satisfying $\Delta_j\in[0,\rho_j]$ for all $j\in[m]$ and $\Delta = \sum_{j=1}^m\Delta_j \leq \mu\sqrt{\dim\mK}\Vert{\w_i-\w^*}\Vert$ such that, $\w^*$ can be written as $\w^* = \sum_{j=1}^m(\rho_j-\Delta_j)\v_j + \Delta\z$, for some $\z\in\mK$.

Additionally, note that, by defining the point $\p_i = \sum_{j=1}^m(\rho_j-\Delta_j)\v_j + \Delta\u_i$ (i.e., replacing the point $\z$ in the representation above of $\w^*$ with the point $\u_i$ computed by the LOO call on iteration $i$ of the algorithm), we have that $\langle{\p_i-\w_i, \nabla{}\Phi_t(\w_i)}\rangle \leq \langle{\w^*-\w_i, \nabla\Phi_t(\w_i)}\rangle \leq -g_i$, where the last inequality is due to convexity of $\Phi_t(\cdot)$.
Thus, we have that
\begin{align}\label{eq:afwproof:1}
-g_i &\geq \langle{\p_i-\w_i, \nabla{}\Phi_t(\w_i)}\rangle= \sum_{j=1}^m\Delta_j\langle{\u_i-\v_j, \nabla{}\Phi_t(\w_i)}\rangle \nonumber \\
&\geq \sum_{j=1}^m\Delta_j\langle{\u_i-\v_1, \nabla{}\Phi_t(\w_i)}\rangle \nonumber \\
&= \Delta\langle{\u_i-\w_i, \nabla{}\Phi_t(\w_i)}\rangle +  \Delta\langle{\w_i-\v_1, \nabla{}\Phi_t(\w_i)}\rangle \nonumber \\
& \geq 2\Delta\langle{\s_i, \nabla\Phi_t(\w_i)}\rangle.
\end{align}


It follows that for any $\zeta > 0$ such that $\zeta\Delta \leq 1$, and when either the Frank-Wolfe direction was chosen or the away direction was chosen  with $\gamma_i < \gamma_{\max}$ (i.e., not a drop step),
\begin{align}\label{eq:afwproof:2}
\Phi_t(\w_{i+1})  &\underset{(a)}{\leq} {\min}_{\gamma\in[0,1]}\Phi_t(\w_i + \gamma\s_i) \leq\Phi_t(\w_i + \zeta\Delta\s_i) \nonumber \\
&\underset{(b)}{\leq}\Phi_t(\w_i) + \zeta\Delta\langle{\s_i, \nabla{}\Phi_t(\w_i)}\rangle + \frac{\zeta^2\Delta^2\beta_t\Vert{\s_i}\Vert^2}{2} \nonumber \\
&\underset{(c)}{\leq} \Phi_t(\w_i) - \frac{\zeta}{2}g_i + \frac{\zeta^2\beta_tD^2}{2}\left({\mu\sqrt{\dim\mK}\Vert{\w_i-\w^*}\Vert}\right)^2 \nonumber \\
&\underset{(d)}{\leq} \Phi_t(\w_i) - \frac{\zeta}{2}g_i + \zeta^2\mu^2D^2\dim\mK{}\cdot{}g_i,
\end{align}
where (a) follows from the use of line-search, and the convexity of $\Phi_t(\cdot)$ which implies that in the case of an away-step with  $\gamma_i < \gamma_{\max}$, $\gamma_i$ is in particular the minimizer of $\Phi_t(\w_i+\gamma\s_i)$ over the entire $\reals$, (b) follows from the smoothness of $\Phi_t(\cdot)$, (c) follows from Eq. \eqref{eq:afwproof:1} and plugging-in the upper-bound on $\Delta$ listed above, and (d) follows from the $\beta_t$-strong convexity of $\Phi_t(\cdot)$.
 
Denoting $\kappa = \mu^2D^2\dim\mK$, we have that for $\zeta = \min\{1, 1/(4\kappa)\}$ (note $\Delta \leq 1$, and thus this indeed satisfies $\zeta\Delta \leq 1$), by subtracting $\Phi_t(\w^*)$ from both sides of \eqref{eq:afwproof:2} we get that,
\begin{align*}
g_{i+1} \leq \left({1 - \min\left\{\frac{1}{4},\frac{1}{16\kappa}\right\}}\right)g_i.
\end{align*}

Now, a generic conversion argument from $g_i$ to $q_i$ (see for instance  Theorem 2 in \cite{lacoste2015global}) yields that 
\begin{align}\label{eq:afwproof:3}
g_i \leq \beta_tD^2/2 \quad \Longrightarrow \quad q_{i} \leq D\sqrt{2\beta_tg_i}.
\end{align}
Note that \eqref{eq:afw:1} (with $i=0$, $\gamma=1$, and denoting $\w_0 = \x_{t-1}$) implies that $g_1 \leq \frac{\beta_t{}D^2}{2}$, and thus the bound on $q_i$ in \eqref{eq:afwproof:3} holds for all $i\geq 1$.

Thus, in order to obtain  the second term inside the $\min$ in Eq. \eqref{eq:afw:conv}, we are interested in the number of iteration $N$ until  $g_{N+1} \leq \frac{\nu_t^2}{2\beta_tD^2} = \frac{\lambda_t^2\nu_t^2}{2\beta{}D^2}$. As before, since for any number of iterations $\tau$ the number of drop steps after $\tau$ iterations cannot exceed $(\tau+1)/2$, we have that 
\begin{align*}
N = O\left({\max\left\{1, \mu^2D^2\dim\mK\right\}\log_+\left({\frac{\beta{}D^2}{\lambda_t^2\nu_t}}\right)}\right).
\end{align*}

The second part of the theorem follows by repeating the above arguments
with respect to the polytope $\mF$. Indeed, since $\w^*\in\mF$, $\w_1\in\mF$ and
$\u_i\in\mF$ for all $i\geq 1$, all iterates and all vertices in their
maintained convex decompositions remain in $\mF$. Thus, the above
arguments apply with $D,\mu,\dim\mK$ replaced with
$D_{\mF},\mu_{\mF},\dim\mF$, respectively.

There is only a minor difference in the initialization of the linear rate
argument: since $\x_{t-1}$ need not be in $\mF$, the above upper-bound on $g_1$ need not translate to $\mF$ (i.e., replacing $D$ with $D_{\mF}$ in the bound $g_1 \leq \beta_t{}D^2/2$). This is fixed as follows. If the algorithm does not terminate at $\w_1$, then its first
step is necessarily a Frank--Wolfe step, since the active set consists
only of $\w_1$. By the use of line-search, the smoothness of $\Phi_t$, and the fact that $\u_1$ is also a linear
minimizer over $\mF$, similarly to  the use of \eqref{eq:afw:1} to upper-bound $g_1$ in the previous case, we have that
$g_2 \leq \frac{\beta_t D_{\mF}^2}{2}$.
Hence the linear rate argument above applies from iteration $i=2$
onward, which changes the iteration bound by at most one.
\end{proof}

\bibliography{bibs.bib}
\bibliographystyle{plain}

\end{document}